\documentclass[12pt]{amsart}
\usepackage{amsmath,amssymb,hyperref,geometry,graphicx,amsthm,tikz,enumitem,mathtools,textcomp}
\usepackage[all,cmtip]{xy}

\newtheorem{theorem}[equation]{Theorem}
\newtheorem{lemma}[equation]{Lemma}
\newtheorem{corollary}[equation]{Corollary}
\newtheorem{proposition}[equation]{Proposition}
\theoremstyle{definition}
\newtheorem{example}[equation]{Example}
\newtheorem{remark}[equation]{Remark}
\newtheorem{definition}[equation]{Definition}

\newtheorem{motivatingquestion}[equation]{Motivating Question}

\newtheorem{innercustomthm}{Theorem}
\newenvironment{customthm}[1]
  {\renewcommand\theinnercustomthm{#1}\innercustomthm}
  {\endinnercustomthm}

\numberwithin{equation}{section}

\newcommand{\Z}{\mathbb{Z}}
\newcommand{\R}{\mathbb{R}}
\newcommand{\Q}{\mathbb{Q}}
\newcommand{\C}{\mathbb{C}}
\newcommand{\N}{\mathbb{N}}
\newcommand{\T}{\mathbb{T}}

\newcommand{\PL}{\mathrm{PL}}
\newcommand{\uPL}{\underline{\mathrm{PL}}}
\newcommand{\PE}{\mathrm{PE}}
\newcommand{\uPE}{\underline{\mathrm{PE}}}
\newcommand{\PP}{\mathrm{PP}}
\newcommand{\uPP}{\underline{\mathrm{PP}}}
\renewcommand{\L}{\mathrm{L}}
\newcommand{\Pic}{\uPL}
\renewcommand{\O}{\mathcal{O}}
\renewcommand{\phi}{\varphi}

\newcommand{\Fun}{\mathrm{Fun}}

\newcommand{\J}{\mathcal{J}}

\newcommand{\K}{\mathrm{K}}
\newcommand{\Span}{\mathrm{Span}}
\newcommand{\col}{\colon}
\newcommand{\RR}{\R}
\newcommand{\ZZ}{\Z}

\newcommand{\M}{\mathcal{M}}

\newcommand{\e}{\mathbf{e}}

\newcommand{\Chi}{\raisebox{2bp}{$\chi$}}

\title{Piecewise-exponential functions and Ehrhart fans}

\author[M.~Chan]{Melody Chan}
\address{Department of Mathematics, Brown University, United States of America}
\email{melody\_chan@brown.edu}

\author[E.~Clader]{Emily Clader}
\address{Department of Mathematics, San Francisco State University, United States of America}
\email{eclader@sfsu.edu}

\author[C.~Klivans]{Caroline Klivans}
\address{Department of Mathematics, Brown University, United States of America}
\email{caroline\_klivans@brown.edu}

\author[D.~Ross]{Dustin Ross}
\address{Department of Mathematics, San Francisco State University, United States of America}
\email{rossd@sfsu.edu}

\begin{document}

\begin{abstract}
This paper studies rings of integral piecewise-exponential functions on rational fans.  Motivated by lattice-point counting in polytopes, we introduce a special class of unimodular fans called Ehrhart fans, whose rings of integral piecewise-exponential functions admit a canonical linear functional that behaves like a lattice-point count. In particular, we verify that all complete unimodular fans are Ehrhart and that the Ehrhart functional agrees with lattice-point counting in corresponding polytopes, which can otherwise be interpreted as holomorphic Euler characteristics of vector bundles on smooth toric varieties.  We also prove that all Bergman fans of matroids are Ehrhart and that the Ehrhart functional in this case agrees with the Euler characteristic of matroids, introduced recently by Larson, Li, Payne, and Proudfoot. A key property that we prove about the Ehrharticity of fans is that it only depends on the support of the fan, not on the fan structure, thus providing a uniform framework for studying $\K$-rings and Euler characteristics of complete fans and Bergman fans simultaneously.
\end{abstract}

\maketitle

\section{Introduction}

The overarching goal of this paper is to develop a unifying combinatorial framework to make sense of two contexts in which the notion of an ``Euler characteristic'' appears.  First, if $\Sigma$ is a unimodular complete fan, so that the associated toric variety $Z_\Sigma$ is smooth and proper, then the holomorphic Euler characteristic is a function
\[
\Chi\colon K(Z_\Sigma) \rightarrow \Z
\]
defined on the Grothendieck $\K$-ring of vector bundles on $Z_{\Sigma}$; combinatorially, $\Chi$ can be interpreted in terms of lattice-point counts in certain polytopes, as we recall below.  The Euler characteristic is not generally well-defined for incomplete fans, but if $\M$ is a matroid, then although its Bergman fan $\Sigma_\M$ is typically not complete, Larson--Li--Payne--Proudfoot \cite{LLPP} recently showed that the $\K$-ring $K(\M)= K(Z_{\Sigma_\M})$ nevertheless admits a function 
\[
\Chi_\M\colon K(\M) \rightarrow \Z
\]
that deserves to be called the Euler characteristic of $\M$.  This can be interpreted geometrically, when $\M$ is realizable, by recognizing that the (proper) wonderful variety of $\M$ embeds into the (nonproper) toric variety $Z_{\Sigma_\M}$, inducing an isomorphism of $K$-rings, and allowing one to define the Euler characteristic on $K(\M)$ using the Euler characteristic of the wonderful variety. Abstracting this description, Larson--Li--Payne--Proudfoot give a purely combinatorial formula for $\Chi_\M$ that applies in the nonrepresentable case, as well.

Put more precisely, then, our goal is to introduce a class of fans---including complete unimodular fans as well as Bergman fans of matroids as special cases---that admit an integer-valued function on their $\K$-ring specializing to the two notions of Euler characteristic in these special cases.

The first key step toward this end is to recall (from work of Brion--Vergne \cite{BrionVergne} and Merkurjev \cite{Merkurjev}) an alternative interpretation of the $\K$-ring of a smooth toric variety in terms of piecewise-exponential functions.  Roughly, for a unimodular fan $\Sigma$ in $N_\R$ with respect to a fixed lattice $N$, there is an additive group $\PL(\Sigma)$ of integral piecewise-linear functions on $\Sigma$, and the $\Z$-algebra of $\PE(\Sigma)$ of piecewise-exponential functions is generated by functions of the form $e^f$ for $f \in \PL(\Sigma)$.  If $\uPL(\Sigma)$ denotes the quotient of $\PL(\Sigma)$ by the subgroup of integral linear functions, and $\uPE(\Sigma)$ denotes the corresponding quotient of $\PE(\Sigma)$ by the relations $e^\ell = 1$ for integral linear functions $\ell$, then the results of \cite{BrionVergne, Merkurjev} show that there is a natural isomorphism
\[
K(Z_\Sigma) \cong \uPE(\Sigma)
\]
for any unimodular fan $\Sigma$.  Moreover, the map $[f] \mapsto [e^f]$ is an injection of $\uPL(\Sigma)$ to the group of units in $\uPE(\Sigma)$; in algebro-geometric terms, it corresponds under the above isomorphism to the Picard group of $\Sigma$.

\begin{figure}[t]
\begin{tikzpicture}[scale=1.5]
\draw[thick,fill=green!20, fill opacity=.5] (-1,-1) -- (1,-1) -- (1,0) -- (0,1) -- (-1,1) -- (-1,-1);

\draw[ultra thick, gray, opacity=0.5] (0,0) -- (0,1.25);
\draw[ultra thick, gray, opacity=0.5] (0,0) -- (0.75, 0.75);
\draw[ultra thick, gray, opacity=0.5] (0,0) -- (1.25, 0);
\draw[ultra thick, gray, opacity=0.5] (0,0) -- (0,-1.25);
\draw[ultra thick, gray, opacity=0.5] (0,0) -- (-1.25,0);

\node at (0,0) {$\bullet$};
\node at (1,0) {$\bullet$};
\node at (1,-1) {$\bullet$};
\node at (0,1) {$\bullet$};
\node at (-1,1) {$\bullet$};
\node at (0,-1) {$\bullet$};
\node at (-1,-1) {$\bullet$};
\node at (-1,0) {$\bullet$};
\node[] at (-.12,-.12) {$0$};
\end{tikzpicture}
\hspace{50bp}
\begin{tikzpicture}[scale=1.5]
\draw[thick,fill=green!20, fill opacity=.5] (-1,-1) -- (1,-1) -- (-1,1) -- (-1,-1);

\draw[ultra thick, gray, opacity=0.5] (0,0) -- (0,1.25);
\draw[ultra thick, gray, opacity=0.5] (0,0) -- (0.75, 0.75);
\draw[ultra thick, gray, opacity=0.5] (0,0) -- (1.25, 0);
\draw[ultra thick, gray, opacity=0.5] (0,0) -- (0,-1.25);
\draw[ultra thick, gray, opacity=0.5] (0,0) -- (-1.25,0);

%\draw[thick] (1,0) -- (0,1);
\node at (0,0) {$\bullet$};
\node at (1,0) {$\bullet$};
\node at (1,-1) {$\bullet$};
\node at (0,1) {$\bullet$};
\node at (-1,1) {$\bullet$};
\node at (0,-1) {$\bullet$};
\node at (-1,-1) {$\bullet$};
\node at (-1,0) {$\bullet$};
\node[] at (-.12,-.12) {$0$};
\end{tikzpicture}
\caption{On the left, a two-dimensional fan $\Sigma$ (with rays drawn in light gray), and the polytope $P_f$ associated to the function $f \in \PL(\Sigma)$ taking value $1$ at each ray generator.  On the right, the polytope $P_{f- \delta_\rho}$, where $\rho$ is the diagonal ray of $\Sigma$. The two lattice points appearing in the left polytope but not in the right polytope are precisely the lattice points of the face $P_f^\rho$.}
\label{fig:intro}
\end{figure}
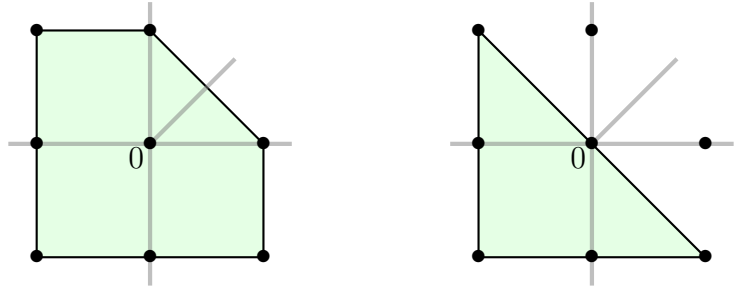

In the case where $\Sigma$ is both unimodular and complete, the fact that the elements $[e^f]$ additively generate $\uPE(\Sigma)$ implies---in light of the above injection---that the Euler characteristic is determined by a function
\[
\Chi_\Sigma: \uPL(\Sigma) \rightarrow \Z,
\]
and this function has a well-known interpretation in terms of lattice-point counts in polytopes.  Namely, any convex function $f \in \PL(\Sigma)$ induces a lattice polytope
\[
P_f = \{m \in M_\R \; | \; m(u_\rho) \leq f(u_\rho) \; \forall \rho \in \Sigma(1)\},
\]
where $M$ denotes the dual lattice of $N$ and $u_\rho$ denotes the ray generator of $\rho$.  In terms of this polytope, the value of $\Chi_\Sigma$ on $[f] \in \uPL(\Sigma)$ is the number of lattice points in $P_f$:
\[
\Chi_\Sigma([f]) = |P_f \cap M|.
\]
Such counts satisfy a fundamental recursion, which we refer to as the Ehrhart recursion: if $P_f^\rho$ is the face of $P_f$ in the direction of some ray $\rho$, and if $\delta_\rho \in \PL(\Sigma)$ is the unique piecewise-linear function taking value $1$ on the ray generator $u_\rho$ and value $0$ at all other ray generators, then
\[
|P_f \cap M| = |P_{f- \delta_\rho} \cap M | + |P_f^\rho \cap M|,
\]
so long as $f - \delta_\rho$ is still convex.  As illustrated in Figure~\ref{fig:intro}, this equation simply reflects the fact that the polytope $P_{f - \delta_\rho}$ is obtained from $P_f$ by ``sliding'' the face $P_f^\delta$ inward by one lattice step, which removes only the lattice points in $P_f^\rho$ from $P_f$.  Phrased fan-theoretically, if $\Sigma^\rho$ denotes the star fan of $\Sigma$ at $\rho$, then there is an induced element $[f]^{\rho} \in \uPL(\Sigma^\rho)$, and the Ehrhart recursion becomes
\begin{equation}
    \label{eq:Ehrhartrecursion}
    \Chi_\Sigma([f]) = \Chi_\Sigma([f - \delta_\rho]) + \Chi_{\Sigma^\rho}([f]^{\rho}).
\end{equation}

Motivated by this setting, we recursively define an {\bf Ehrhart fan} to be a unimodular fan admitting a function $\Chi_\Sigma\colon \uPL(\Sigma) \rightarrow \Z$ satisfying \eqref{eq:Ehrhartrecursion} for all $[f] \in \uPL(\Sigma)$ and $\rho\in\Sigma(1)$, normalized by $\Chi_\Sigma(0)=1$. 
 We make this notion precise in Definition~\ref{def:Ehrhart}, and then we prove that the function $\Chi_\Sigma$ on an Ehrhart fan is unique and polynomial.  The following two theorems---stated in more detail in the body of the paper---show that Ehrhart fans indeed include complete unimodular fans and Bergman fans as special cases:
 
\begin{customthm}{A}(Theorem~\ref{thm:complete}) \label{thm:A} Let $\Sigma$ be a complete unimodular fan.  Then $\Sigma$ is Ehrhart, and when $f \in \PL(\Sigma)$ is convex, $\Chi_\Sigma([f])$ is given by the lattice-point count:
\[
\Chi_\Sigma([f])= |P_f \cap M|.
\]
\end{customthm}

\begin{customthm}{B}(Theorem~\ref{thm:matroid}) \label{thm:B} Let $\Sigma_\M$ be the Bergman fan of a matroid $\M$.  Then $\Sigma_\M$ is Ehrhart, and under an appropriate inclusion $\uPL(\Sigma_\M) \subseteq K(\M)$, the Ehrhart function $\Chi_{\Sigma_\M}$ agrees with the Euler characteristic $\Chi_\M$ of Larson--Li--Payne--Proudfoot:
\[
\Chi_{\Sigma_\M}([f]) = \Chi_\M([f]).
\]
\end{customthm}

These theorems hint at two related features of Ehrharticity that are not at all obvious.  First, since lattice-point counts do not depend on the particular complete fan at hand, only on the piecewise-linear functions, Theorem~\ref{thm:A} suggests that Ehrharticity depends only on the support of $\Sigma$, not on $\Sigma$ itself.  This is indeed the case:

\begin{customthm}{C}(Theorem~\ref{thm:support})
\label{thm:C} Let $\Sigma_1$ and $\Sigma_2$ be unimodular fans with $|\Sigma_1| = |\Sigma_2|$.  Then $\Sigma_1$ is Ehrhart if and only if $\Sigma_2$ is Ehrhart.  Moreover, if $f \in \PL(\Sigma_1) \cap \PL(\Sigma_2)$, then
\[
\Chi_{\Sigma_1}([f]) = \Chi_{\Sigma_2}([f]).
\]
\end{customthm}

Second, since Euler characteristics of complete fans and Bergman fans are well-defined on the K-rings of the associated toric varieties, not just Picard groups, Theorems~\ref{thm:A} and \ref{thm:B} suggest that the Ehrhart function of an Ehrhart fan $\Sigma$ should induce a well-defined function
\[\widetilde{\Chi}_\Sigma: \uPE(\Sigma) \rightarrow \Z\]
given by $\widetilde{\Chi}_\Sigma([e^f]) = \Chi_\Sigma([f])$.  This is also the case, and in fact, Theorem~\ref{thm:C} is the key to its proof.  To understand this, it is helpful to note---as we discuss in Section \ref{subsec:uPE(X)presentation}---that the ring $\uPE(\Sigma)$ has relations 
\begin{equation}
\label{eq:introrelation}
[e^f] + [e^g] = [e^{\text{max}(f,g)}] + [e^{\text{min}(f,g)}]
\end{equation}
for any $f,g \in \PL(\Sigma)$.  However, even if $f$ and $g$ are piecewise-linear on $\Sigma$, their maximum and minimum are generally only piecewise-linear on a refinement of $\Sigma$.  In light of this, the relation \eqref{eq:introrelation} is most naturally viewed as a relation in the ring $\uPE(X)$ of piecewise-exponential functions on {\it any} unimodular fan with support $X= |\Sigma|$ (that is, the direct limit of the rings $\uPE(\Sigma')$ for $|\Sigma'| = X$); in fact, as we discuss in Theorem~\ref{thm:K(X)presentation}, relations of this form are, in a precise sense, the only relations in $\uPE(X)$.  Theorem~\ref{thm:C} implies that there is a well-defined Ehrhart function
\[\Chi_X: \uPL(X) \rightarrow \Z\]
on the analogously-defined ring of piecewise-linear functions on $X$, and it is not much additional work to see that this function respects the relation \eqref{eq:introrelation}, so one finds (Theorem~\ref{thm:chiK}) that there is a well-defined function
\[\widetilde{\Chi}_X: \uPE(X) \rightarrow \Z\]
given by $\widetilde{\Chi}_X([e^f]) = \Chi_X([f])$ for any $f \in \PL(X)$.  

In the case where $\Sigma$ is complete and thus $X= |\Sigma| = N_\R$, work of Brion \cite{Brion} and Pukhlikov--Khovanskii \cite{PukhlikovKhovanskii} shows that
\[\uPE(X) \cong \Pi(M),\]
where $\Pi(M)$ denotes the polytope algebra (studied by McMullen \cite{McMullen} and others), which is the free abelian group generated by lattice polytopes $P \subseteq M_\R$ up to translation, endowed with a ring structure via Minkowski sum.  Under this isomorphism, the function
\[\widetilde{\Chi}_{X}: \Pi(M) \rightarrow \Z\]
is the (well-defined) lattice-point count on the polytope algebra; furthermore, by the discussion above, upon restricting to piecewise-exponential functions on a fixed complete fan $\Sigma$, this agrees with the holomorphic Euler characteristic of $Z_\Sigma$.  In the case where $\Sigma$ is the Bergman fan of a matroid $\M$, on the other hand, Theorem~\ref{thm:B} shows that
\[\widetilde{\Chi}_\Sigma: K(\M) \rightarrow \Z\]
is Larson--Li--Payne--Proudfoot's Euler characteristic.  Thus, Ehrhart fans indeed achieve the goal of unifying the two notions of Euler characteristic described at the outset.

\subsection*{Plan of the paper}

The paper begins with a substantial amount of background on rings of piecewise-exponential functions; some of this material is likely familiar to experts, but because existing treatments of these rings are in a different context from ours, we aim to build up all of the necessary properties in a self-contained way.  In particular, Section~\ref{sec:PE} gives the definitions of the rings of piecewise-linear and piecewise-exponential functions both on fans and on fan-supports, and describes an explicit presentation of $\uPE(X)$ in terms of the relations \eqref{eq:introrelation}.  In Section~\ref{sec:Ehrhart}, we define Ehrhart fans, and we prove some of the first key properties: polynomiality of the Ehrhart function, Ehrhart fans are balanced (though not all balanced fans are Ehrhart), and an interpretation of the leading term of the Ehrhart function as a ``volume polynomial'' of the fan.  Section~\ref{sec:complete} gives the proof that complete fans are Ehrhart, Section~\ref{sec:independence} that Ehrharticity depends only on the support of a fan, and finally, Section~\ref{sec:matroidfans} gives the proof that Bergman fans of matroids are Ehrhart.

\subsection*{Acknowledgments}  

This work began while E.C. and D.R. were visiting Brown University for sabbatical; they warmly acknowledge the Brown Department of Mathematics for their support and for providing a welcoming and stimulating research environment. The authors thank Matt Larson for helpful feedback on an early draft. D.R. also thanks Jupiter Davis and Ian Wallace for enlightening conversations related to what ultimately became the notion of an Ehrhart fan.  M.C.~was supported by NSF CAREER DMS--1844768, FRG DMS--2053221, and DMS--2401282. E.C. was supported by NSF CAREER DMS--2137060.  D.R. was supported by NSF DMS--2302024 and DMS--2001439.

\section{Piecewise-exponential functions}
\label{sec:PE}

In this section, we introduce the primary algebraic objects of study in this paper: the groups of integral piecewise-linear functions and the rings of integral piecewise-exponential functions on rational fans and their supports.  These algebraic objects have previously been studied in various different guises and closely-related contexts; for example, they were studied explicitly in Brion's work \cite{Brion} and implicitly in the work of Khovanski\u{\i} and Pukhlikov \cite{PukhlikovKhovanskii}. We begin with a brief discussion of our conventions regarding rational fans.

\subsection{Fan conventions}
\label{subsection:conventions}

Throughout this paper, $N$ denotes a free abelian group of rank $n$, and $N_\R=N\otimes_\Z\R$ denotes the associated vector space. We let $M=\mathrm{Hom}_\Z(N,\Z)$ denote the dual of $N$ and $M_\R=\mathrm{Hom}_\R(N_\R,\R)$ the dual of $N_\R$. Note that $M_\R$ is the space of linear functions on $N_\R$, and $M\subseteq M_\R$ is the subset of linear functions that take integer values on all elements of the lattice $N\subseteq N_\R$.  There is a natural pairing $M_\R \times N_\R \rightarrow \R$ that we denote by $\langle v, u \rangle \coloneq v(u)$.

Closed and open {\bf rational halfspaces} in $N_\R$ are subsets of the form
\[
H_m=\{v\in N_\R\mid \ \langle m,v\rangle\geq 0\}\;\;\;\text{ and }\;\;\;H_m^\circ=\{v\in N_\R\mid \ \langle m,v\rangle> 0\}
\]
for some $m\in M$. A {\bf rational cone} is a finite intersection of closed rational halfspaces. The {\bf dimension} of a rational cone is the dimension of the smallest linear subspace of $N_\R$ containing it. If $\sigma$ is a rational cone, then a subset $\tau\subseteq\sigma$ is a {\bf face} of $\sigma$, written $\tau\preceq\sigma$, if 
\[
\tau=\sigma\cap H_m\;\; \text{ for some $m\in M$ with }\;\;\sigma\cap H_m^\circ=\emptyset.
\]
A {\bf rational fan} in $N_\R$ is a finite collection $\Sigma$ of rational cones in $N_\R$ such that
\begin{enumerate}
    \item If $\sigma\in\Sigma$ and $\tau\preceq \sigma$, then $\tau\in\Sigma$, and
    \item If $\sigma_1,\sigma_2\in\Sigma$, then $\sigma_1\cap\sigma_2\preceq\sigma_1$ and $\sigma_1\cap\sigma_2\preceq\sigma_2$.
\end{enumerate}
A {\bf subfan} is a subset of a fan that is, itself, a fan. The {\bf dimension} of a fan is the largest dimension of any cone in the fan, and a fan is said to be {\bf pure} if every inclusion-maximal cone has the same dimension. Given a rational fan $\Sigma$, we let $\Sigma(k)$ denote the set of $k$-dimensional cones of $\Sigma$.  The elements of $\Sigma(1)$ are called {\bf rays}.

A rational fan is {\bf pointed} if it contains the origin (the unique zero-dimensional rational cone). All one-dimensional cones in a pointed rational fan are rays emanating from the origin, and for each such ray $\rho\in\Sigma$, the {\bf ray generator} $u_\rho\in N$ is defined as the unique element in $\rho \cap N$ that generates the latter as a monoid.  A pointed fan is said to be {\bf simplicial} if the number of rays in each cone is equal to the dimension of the cone, and it is said to be {\bf unimodular} if the set of ray generators of each cone can be extended to a $\Z$-basis of $N$. The {\bf support} of a rational fan $\Sigma$, denoted $|\Sigma|$, is the union of all cones $\sigma\in\Sigma$, and we say that $\Sigma$ is {\bf complete} if $|\Sigma|=N_\R$.

If $\Sigma$ and $\Sigma'$ are two fans with the same support, we say that $\Sigma'$ {\bf refines} $\Sigma$ if every cone of $\Sigma'$ is contained in some cone of $\Sigma$.    It is a well-known fact that every rational fan admits a unimodular refinement \cite[Theorem~11.1.9]{CoxLittleSchenck}. Furthermore, any two rational fans $\Sigma_1$ and $\Sigma_2$ naturally admit a common rational refinement
\[
\Sigma'=\{\sigma_1\cap\sigma_2\mid \sigma_1\in\Sigma_1,\;\sigma_2\in\Sigma_2\},
\]
which can then be refined into a common unimodular refinement.

Given a rational fan $\Sigma$ and a cone $\sigma\in\Sigma$, we define the {\bf neighborhood} of $\sigma$ in $\Sigma$ as the subfan $\Sigma^{(\sigma)} \subseteq \Sigma$ comprised of the cones containing $\sigma$, and all of the faces of those cones; that is,
\[
\Sigma^{(\sigma)}\coloneq\{\tau\in\Sigma\mid \tau\preceq \sigma'\text{ for some }\sigma'\in\Sigma\text{ with }\sigma\preceq\sigma'\}.
\]
Next, consider the free abelian group %(free) quotient group
\[
N^\sigma\coloneq N/\Span_\Z(\sigma\cap N)
\]
with corresponding vector space $N_\R^\sigma:=N^\sigma\otimes_\Z\R$.  The {\bf star fan} of $\Sigma$ at $\sigma$ is defined as the image of the neighborhood $\Sigma^{(\sigma)}$ under the quotient map $N_\R\rightarrow N^\sigma_\R$:
\[
\Sigma^\sigma\coloneq\{\overline\tau\mid\tau\in\Sigma^{(\sigma)}\},
\]
where we write $\overline\tau$ to denote the image of $\tau\subseteq N_\R$ under the quotient map $N_\R \rightarrow N_\R^\sigma$.  If the original fan $\Sigma$ is unimodular, then $\Sigma^{\sigma}$ is also unimodular with respect to the lattice $N^\sigma\subseteq N_\R^\sigma$.

\subsection{Piecewise-linear functions}

The primary algebraic objects studied in this paper are built out of integral piecewise-linear functions on rational fans. To define them, we view an element $m \in M$ as an integral linear function on $N_\R$---that is, a linear function $N_\R \rightarrow \R$ that takes $N$ to $\Z$.  

\begin{definition}
\label{def:PL}
Let $\Sigma$ be a rational fan in $N_\R$. A function $f\colon |\Sigma|\rightarrow\R$ is \textbf{integral piecewise-linear on $\Sigma$} if, for every cone $\sigma\in\Sigma$, there exists an integral linear function $m\in M$ such that
\[
f|_\sigma=m|_\sigma.
\]
Denote by $\PL(\Sigma)$ the group of integral piecewise-linear functions on $\Sigma$.
\end{definition}

Integral piecewise-linear functions on $\Sigma$ are uniquely determined by their values at the ray generators of any pointed fan $\Sigma'$ that refines $\Sigma$.  In particular, this implies that $\PL(\Sigma)$ is a subgroup of $\Z^{\Sigma'(1)}$, from which it follows that $\PL(\Sigma)$ is a finitely-generated free abelian group. However, even when $\Sigma$ is simplicial, not every choice of integer values at ray generators corresponds to a piecewise-linear function, as the next example illustrates.

\begin{example}\label{example:notlinear}
Let $N=\ZZ^2$, and let $\Sigma$ be the two-dimensional fan in $N_\RR \cong \RR^2$ with one maximal cone $\sigma$ spanned by the two rays $v_1=(1,0)$ and $v_2=(1,2)$. There does not exist an integral piecewise-linear function $f\in\PL(\Sigma)$ such that $f(v_1)=0$ and $f(v_2)=1$. Indeed, if such a function existed, then linearity on $\sigma$ would imply that its value at
\[
(1,1)=\frac{1}{2}(1,0)+\frac{1}{2}(1,2) \in \sigma
\]
is $1/2$, but this is not an integer.
\end{example}

In light of this example, we often restrict in what follows to unimodular fans.  In particular, if $\Sigma$ is unimodular, then specifying an integral piecewise-linear function on $\Sigma$ is equivalent to choosing an integer value at each ray generator, so
\[
\PL(\Sigma)=\Z^{\Sigma(1)}.
\]

We will primarily be interested in the quotient of $\PL(\Sigma)$ under translation by linear functions, defined as follows. 

\begin{definition}\label{def:uPL}
Let $\Sigma$ be a rational fan in $N_\R$. An element $f\in\PL(\Sigma)$ is \textbf{integral linear on $\Sigma$} if there exists $m\in M$ such that
\[
f=m|_{|\Sigma|}.
\]
Let $\L(\Sigma)\subseteq\PL(\Sigma)$ denote the subgroup of integral linear functions on $\Sigma$, and denote by
\[
\uPL(\Sigma)\coloneq\frac{\PL(\Sigma)}{\L(\Sigma)}
\]
the \textbf{group of integral piecewise-linear functions on $\Sigma$ up to translation}.
\end{definition}

We note that $\uPL(\Sigma)$ is often denoted by $A^1(\Sigma)$ in the literature---as in \cite{AHK}, for example---and it can naturally be identified with the class group of the corresponding toric variety $Z_\Sigma$, which is equal to the Picard group when $\Sigma$ is unimodular \cite[Chapter 4]{CoxLittleSchenck}.

The abelian group $\uPL(\Sigma)$ is finitely-generated, but it is generally not free, even for unimodular fans; the next example illustrates this phenomenon explicitly.

\begin{example}\label{example:torsion}
Let $N=\ZZ^2$, and let $\Sigma$ be the one-dimensional fan in $N_\RR\cong \R^2$ with two rays, one spanned by $v_1=(1,0)$ and one spanned by $v_2=(1,2)$. There is a unique piecewise-linear function $f$ such that $f(v_1)=0$ and $f(v_2)=1$. This piecewise-linear function is not the restriction of an integral linear function on $\R^2$ (as we saw in Example~\ref{example:notlinear}), but $2f$ \emph{is} the restriction of an integral linear function; specifically, $2f$ is the restriction of the integral linear function $(x,y)\mapsto y$. In other words,
\[
[f]\neq 0\;\;\;\text{ but }\;\;\;2[f]=0
\]
in $\uPL(\Sigma)$.  Thus, $\uPL(\Sigma)$ has torsion.
\end{example}

In what follows, we generally drop the adjective ``integral" from the terminology in Definitions~\ref{def:PL} and \ref{def:uPL}, since we only ever consider functions $N_\R \rightarrow \R$ that take integer values on $N$. However, the previous examples highlight that the integrality constraint is important to keep in mind, as it can affect whether a function on $|\Sigma|$ is the restriction of a linear function on $N_\R$.

When working with piecewise-linear functions, it is often useful to work with all fans having the same support simultaneously; for instance, it will be important to consider functions of the form $\max(f,g)$ and $\min(f,g)$ for $f,g \in \PL(\Sigma)$, and these are generally only piecewise-linear on a refinement of $\Sigma$.  This motivates the following definition.

\begin{definition}
Let $X\subseteq N_\R$ be the support of a rational fan. A function $f\colon X\rightarrow\R$ is \textbf{(integral) piecewise-linear on $X$} if there is a rational fan $\Sigma$ with $|\Sigma|=X$ such that $f\in\PL(\Sigma)$. Denote by $\PL(X)$ the group of piecewise-linear functions on $X$ and by $\L(X)\subseteq \PL(X)$ the group of linear functions on $X$. Denote by
\[
\uPL(X)\coloneq \frac{\PL(X)}{\L(X)}
\]
the \textbf{group of (integral) piecewise-linear functions on $X$ up to translation}.
\end{definition}

\begin{remark}
     $\PL(X)$ is closed (and therefore forms a group) under addition of functions.  To see this, suppose that $\Sigma_1$ and $\Sigma_2$ are rational fans with equal support $X$ and $f_1\in \PL(\Sigma_1)$ and $f_2\in \PL(\Sigma_2)$. We may then choose a common refinement $\Sigma$ of $\Sigma_1$ and $\Sigma_2$, and since each cone of $\Sigma$ is contained in a cone of $\Sigma_1$ and a cone of $\Sigma_2$, it follows that $f_1$ and $f_2$ are each linear on the cones $\Sigma$. Thus, the sum $f_1+f_2$ is piecewise-linear on $\Sigma$, and since $|\Sigma|=X$, this implies that $f_1+f_2\in\PL(X)$.
\end{remark}

For any rational fan $\Sigma$ with $|\Sigma|=X$, the natural inclusion $\PL(\Sigma)\subseteq\PL(X)$ gives an identification of subgroups $\L(\Sigma)=\L(X)$. It follows that there is an inclusion of quotient groups
\[
\uPL(\Sigma)\subseteq \uPL(X).
\]
More generally, if $\Sigma'$ is a refinement of $\Sigma$, then there is an inclusion
\[
\uPL(\Sigma)\subseteq \uPL(\Sigma'),
\]
and we can view $\uPL(X)$ as the direct limit
\[
\uPL(X)=\varinjlim \uPL(\Sigma),
\]
where the underlying directed system is the set of rational fans with support $X$, directed by refinement. In fact, since every rational fan admits a unimodular refinement, it follows that $\uPL(X)$ can be viewed as a direct limit over the smaller directed set of unimodular fans with support $X$.

\subsection{Piecewise-exponential functions}
\label{subsec:PE}

We now discuss piecewise-exponential functions, which naturally give rise to a ring containing the group of piecewise-linear functions.

\begin{definition}
Let $\Sigma$ be a rational fan in $N_\R$. A function $F\colon |\Sigma|\rightarrow\R$ is \textbf{(integral) piecewise-exponential on $\Sigma$} if there exists $f\in\PL(\Sigma)$ such that
\[
F=e^f.
\]
Denoting by $\Fun(|\Sigma|,\R)$ the ring of all functions $|\Sigma| \rightarrow \R$, the \textbf{ring of (integral) piecewise-exponential functions on $\Sigma$}, denoted $\PE(\Sigma)$, is the subring of $\Fun(|\Sigma|,\R)$ generated by all (integral) piecewise-exponential functions on $\Sigma$.
\end{definition}

By definition, notice that any element of $\PE(\Sigma)$ can be expressed as a difference 
\[
\sum_{i=1}^k e^{f_i}-\sum_{j=1}^\ell e^{g_j}
\]
for (not necessarily distinct) functions $f_1, \ldots, f_k,g_1,\dots,g_\ell \in \PL(\Sigma)$.  However, this expression is not unique.  The following result characterizes when two such expressions are equal in $\PE(\Sigma)$.

\begin{lemma}\label{lemma:expfunctions}
Let $\Sigma$ be a rational fan. We have equality
\[
\sum_{i=1}^ke^{f_i}=\sum_{j=1}^\ell e^{g_j}\in \PE(\Sigma)
\]
if and only if, for every $v\in |\Sigma|$, there is an equality of multisets
\[
\{f_i(v)\mid i=1,\dots,k\}=\{g_j(v)\mid j=1,\dots,\ell\}.
\]
\end{lemma}

\begin{proof}
The ``if'' direction is immediate. To prove the ``only-if'' direction by contrapositive, suppose that there exists some $v\in|\Sigma|$ such that the two multisets are not equal. Without loss of generality, assume that the elements of the multisets are ordered in nondecreasing order.  By successively removing the largest term from each multiset if they agree, one eventually arrives at two new multisets
\[
\{f_1(v),\dots,f_{k-d}(v)\}\;\;\;\text{ and }\;\;\;\{g_1(v),\dots,g_{\ell-d}(v)\}
\]
such that, without loss of generality, $f_{k-d}(v)>g_{\ell-d}(v)\geq\dots\geq g_1(v)$ and
\begin{align*}
\sum_{i=1}^ke^{f_i(v)}-\sum_{j=1}^\ell e^{g_j(v)}=\sum_{i=1}^{k-d}e^{f_i(v)}-\sum_{j=1}^{\ell-d} e^{g_j(v)}.
\end{align*}
Note that it is possible to have $\ell-d=0$, in which case the final sum above vanishes. Either way, evaluating at $\lambda v$ for $\lambda\gg 0$, the term
\[
e^{f_{k-d}(\lambda v)}=e^{\lambda f_{k-d}(v)}
\]
dominates all other terms, so the difference cannot be zero.
\end{proof}

In light of Lemma~\ref{lemma:expfunctions}, we can view the piecewise-exponential construction simply as a formal method for turning addition in $\PL(\Sigma)$ into multiplication in $\PE(\Sigma)$. In particular, it follows from Lemma~\ref{lemma:expfunctions} that the natural function
\begin{align*}
\PL(\Sigma)&\rightarrow \PE(\Sigma)\\
f&\mapsto e^f,
\end{align*}
is an injective group homomorphism from the additive group $\PL(\Sigma)$ to the multiplicative group of units in $\PE(\Sigma)$.

In the same way that we set $[\ell] = 0 \in \uPL(\Sigma)$ for any linear function $\ell$, we are interested in the quotient of $\PE(\Sigma)$ by the relations $e^\ell = e^0$, as follows.

\begin{definition}
Let $\Sigma$ be a rational fan in $N_\R$. Denote by
\begin{equation}
    \label{eq:uPESigmadef}
\uPE(\Sigma)\coloneq\frac{\PE(\Sigma)}{\big\langle e^\ell-1 \mid \ell\in\L(\Sigma)\big\rangle}
\end{equation}
the \textbf{ring of (integral) piecewise-exponential functions on $\Sigma$ up to translation}.
\end{definition}

By construction, we have $[e^f] = 1 \in \uPE(\Sigma)$ if $f \in \L(\Sigma)$; what is somewhat less obvious is that the converse is true in the following sense.

\begin{lemma}
\label{lem:1inuPE}
Let $f \in \PL(\Sigma)$.  Then $[e^f] = 1 \in \uPE(\Sigma)$ if and only if $f \in \L(\Sigma)$.
\end{lemma}
\begin{proof}
As explained above, the only direction that requires proof is the ``only if'' direction, so suppose that $[e^f] = 1 \in \uPE(\Sigma)$.  This means that $e^f - 1 \in \langle e^\ell - 1 \; | \; \ell \in \L(\Sigma) \rangle$, and unpacking the elements of this ideal, we find that
\[e^f - e^0 = \sum_{i=1}^k \left(e^{f_i + \ell_i} - e^{f_i}\right)\]
for some $f_1, \ldots, f_k \in \PL(\Sigma)$ and $\ell_1, \ldots, \ell_k \in \L(\Sigma)$.  Rearranging and applying Lemma~\ref{lemma:expfunctions}, there is an equality of multisets
\[\{f(v), f_1(v), \ldots, f_k(v)\} = \{0, f_1(v) + \ell_1(v), \ldots, f_k(v) + \ell_k(v)\}\]
for all $v \in |\Sigma|$.  This, in particular, implies that the sums of the elements in these multisets agree, meaning that
\[f(v) = \ell_1(v) + \cdots + \ell_k(v)\]
for all $v \in |\Sigma|$.  That is, $f = \ell_1 + \cdots + \ell_k$, so the lemma is proved.
\end{proof}

It follows from Lemma~\ref{lem:1inuPE} that the natural function
\begin{align*}
\uPL(\Sigma)&\rightarrow \uPE(\Sigma)\\
[f]&\mapsto [e^f],
\end{align*}
is an injective group homomorphism from the additive group $\uPL(\Sigma)$ to the multiplicative group of units in $\uPE(\Sigma)$.

Building on Example~\ref{example:torsion}, we now illustrate how the ring $\uPE(\Sigma)$ may have additive torsion, even for a unimodular fan $\Sigma$.

\begin{example}\label{example:petorsion}
Let $N=\ZZ^2$, and let $\Sigma$ be the one-dimensional fan in $N_\R\cong \R^2$ considered in Example~\ref{example:torsion}. In that example, we found a piecewise-linear function $f\in\PL(\Sigma)\setminus\L(\Sigma)$ such that $2f\in\L(\Sigma)$. Thus, in $\uPE(\Sigma)$, we obtain an element
\[
[e^f]\neq 1\;\;\;\text{ such that }\;\;\;[e^f]^2=1.
\]
(The fact that $[e^f] \neq 1$ follows from Lemma~\ref{lem:1inuPE}, since $f \notin \L(\Sigma)$.)  We now use the function $f$ to show that $\uPE(\Sigma)$ has torsion as an additive group.

For brevity, denote a piecewise-linear function on $\Sigma$ as a pair $(a_1,a_2)\in\Z^2$, where $a_i$ is the value of the function at $v_i$; for instance, the function $f$ from Example~\ref{example:torsion} is represented as the pair $(0,1)$. Notice that a piecewise-linear function represented by $(a_1,a_2)$ is the restriction of a linear function if and only if $a_1-a_2\in2\Z$. 
Within the quotient ring $\uPE(\Sigma)$, we have, for example,
\[[e^f]=[e^{(0,1)}]=[e^{(1,0)}].\]
Using this observation, we compute
\begin{align*}
2[e^f-1]&=[e^{(0,1)}]+[e^{(0,1)}]-[e^{(0,0)}]-[e^{(0,0)}]\\
&=[e^{(0,1)}]+[e^{(1,0)}]-[e^{(0,0)}]-[e^{(0,0)}]\\
&=[e^{(1,1)}]+[e^{(0,0)}]-[e^{(0,0)}]-[e^{(0,0)}]\\
&=[e^{(1,1)}]-[e^{(0,0)}] \\&= 0,
\end{align*}
where the third equality follows from Lemma~\ref{lemma:expfunctions} and the last equality from Lemma~\ref{lem:1inuPE}.  Therefore, $\uPE(\Sigma)$ has nontrivial torsion.
\end{example}

As in the case of piecewise-linear functions, we introduce notation for working with piecewise-exponential functions on all fans with the same support simultaneously.

\begin{definition}
Let $X\subseteq N_\R$ be the support of a rational fan. A function $F\colon X\rightarrow\R$ is \textbf{(integral) piecewise-exponential on $X$} if there is a rational fan $\Sigma$ with $|\Sigma|=X$ such that $F\in\PE(\Sigma)$. Denote by $\PE(X)$ the \textbf{ring of (integral) piecewise-exponential functions on $X$}, and denote by
\begin{equation}
    \label{eq:uPEXdef}
    \uPE(X)\coloneq\frac{\PE(X)}{\big\langle e^\ell-1 \mid \ell\in\L(X)\big\rangle}
\end{equation}
the \textbf{ring of (integral) piecewise-exponential functions on $X$ up to translation}.
\end{definition}

For any rational fan $\Sigma$ with $|\Sigma|=X$, notice that $\PE(\Sigma)\subseteq\PE(X)$ and that the ideal $\langle e^\ell-1 \mid \ell\in\L(\Sigma)\rangle\subseteq\PE(\Sigma)$ is a subset of the ideal $\langle e^\ell-1 \mid \ell\in\L(X)\rangle\subseteq\PE(X)$ via this inclusion. Thus, there is a natural map
\[
\phi_\Sigma:\uPE(\Sigma)\rightarrow\uPE(X).
\]
In fact, there is a natural map $\uPE(\Sigma)\rightarrow\uPE(\Sigma')$ whenever $\Sigma'$ refines $\Sigma$, and we can view $\uPE(X)$ as the direct limit over these maps:
\[
\uPE(X)=\varinjlim \uPE(\Sigma),
\]
where the limit is either over the directed set of all rational fans or over the directed subset of unimodular fans.

\begin{remark}
The ring of not-necessarily-integral piecewise-exponential functions, which we denote $\uPE(X)_\R$, was previously studied by Brion \cite{Brion} in relation to the polytope algebra and Chow rings of toric varieties. The natural inclusion of integral piecewise-exponential functions $\PE(X)\subseteq\PE(X)_\R$ induces a ring homomorphism from $\uPE(X)$ to $\uPE(X)_\R$, but the nontrivial torsion element constructed in Example~\ref{example:petorsion} is in the kernel, so this homomorphism is not generally an inclusion.
\end{remark}

In the next subsection, we give a concrete presentation of the ring $\uPE(X)$.

\subsection{Presentation of $\uPE(X)$}
\label{subsec:uPE(X)presentation}

Let $\Z\{\PL(X)\}$ denote the group algebra of the group $\PL(X)$, whose elements are formal $\Z$-linear combinations of generators $\{f\}$ for each $f \in \PL(X)$, with multiplication defined by
\begin{equation}
\label{eq:groupalgebra}
\{f\} \cdot \{g\} = \{f+g\}.
\end{equation}
It follows from \eqref{eq:groupalgebra} that there is a surjective ring homomorphism
\begin{align}
\label{eq:tau}
 \tau: \Z\{\PL(X)\} &\rightarrow \PE(X),\\
\nonumber \{f\} &\mapsto e^f
\end{align}
(extended $\Z$-linearly), and in this subsection, we study the kernel of $\tau$ in order to prove the following presentation of the ring $\uPE(X)$.

\begin{theorem}
\label{thm:K(X)presentation}
Let $X$ be the support of a rational fan.  There is an isomorphism of rings
\[
\uPE(X) \cong \frac{\Z\{\PL(X)\}}{\Z\bigg\{ \substack{\vspace{5bp}\{f\} + \{g\} - \{\max(f,g)\} - \{\min(f,g)\},\\     \{f+\ell\} - \{f\}}\;\Big|\;\substack{\vspace{5bp}f,g\in\PL(X),\\     \ell\in \L(X)}\bigg\}},
\]
where the notation $\Z\{ \cdot \}$ in the denominator represents the $\Z$-linear span.
\end{theorem}

Two observations should be made about the right-hand side of the theorem.  First, multiplying  either type of generator appearing in the denominator by an element $\{h\}$ with $h\in\PL(X)$ results in another generator, so the $\Z$-linear span in the denominator is actually an ideal and therefore the quotient is a well-defined ring.  Second, the maximum or minimum of two piecewise linear functions on $X$ is, indeed, a piecewise-linear functions on $X$, so expressions like $\{\max(f,g)\}$ and $\{\min(f,g)\}$ make sense in $\Z\{\PL(X)\}$.  This second claim is the content of the following lemma.

\begin{lemma}
Let $X$ be the support of a rational fan. If $f,g\in\PL(X)$, then 
\[
\max(f,g), \; \min(f,g)\in\PL(X).
\]
\end{lemma}

\begin{proof}
Since any two rational fans admit a common rational refinement, we may choose a fan $\Sigma$ with $|\Sigma|=X$ such that $f,g\in\PL(\Sigma).$ For each cone $\sigma\in\Sigma$, the points
\[
\{v\in\sigma\mid f(v)=g(v)\}
\]
are the intersection of $\sigma$ with a rational hyperplane in $\mathrm{span}_\R(\sigma)$. Upon subdividing each cone $\sigma$ along the corresponding hyperplane, we obtain a new rational fan $\Sigma'$ such that $\max(f,g)$ and $\min(f,g)$ are linear when restricted to each cone of $\Sigma'$. It then follows that $\max(f,g),\min(f,g)\in\PL(\Sigma')\subseteq\PL(X)$.
\end{proof}

Viewing the two types of relations in Theorem~\ref{thm:K(X)presentation} separately, relations of the second type arise from the definition of $\uPE(X)$ as a quotient as in \eqref{eq:uPEXdef}.  The crux of the matter in proving Theorem~\ref{thm:K(X)presentation}, then, is the following result.

\begin{proposition}
\label{prop:K(X)presentation}
Let $X$ be the support of a rational fan.  Then
\[
\PE(X) = \Z\{\PL(X)\}/I\]
where
\[I={\Z\big\{ \{f\} + \{g\} - \{\max(f,g)\} - \{\min(f,g)\}\mid f,g\in\PL(X)\big\}}.
\]
\end{proposition}

\begin{proof}
Equivalently, we prove that $I$ is the kernel of the surjective homomorphism $\tau$ defined by \eqref{eq:tau}.  For one inclusion, note that there is an equality of multisets
\[
\{f(v), g(v)\} = \{\max(f,g)(v), \min(f,g)(v)\}  
\]
for all $v \in X$.  Thus,
\[
e^f + e^g = e^{\max(f,g)} + e^{\min(f,g)} \in \PE(X),
\]
from which it follows that $I\subseteq\ker(\tau)$.

To prove the converse, we first prove the following ``reordering of values'' claim in the quotient ring $\Z\{\PL(X)\}/I$: given $f_1,\ldots, f_m \in \PL(X)$, there exist $\widetilde f_1,\ldots, \widetilde f_m\in \PL(X)$ such that $\sum \{f_i\} = \sum \{\widetilde f_i\} \mod I$ and such that
\[\widetilde f_1(v) \le \cdots \le \widetilde f_m(v) \text{ for all } v\in X.\]
This can be proved by induction on the number of $i \in \{1, \ldots, m-1\}$ such that $f_i(v) > f_{i+1}(v)$ for some $v\in X$.  If there are no such $i$, we are done.  Otherwise, let $i$ be as above, and for each $k=1,\ldots,m,$ let
\[\widetilde f_k = \begin{cases} \min (f_i, f_{i+1}) &\text{ if }k=i,\\ \max (f_i, f_{i+1}) & \text{ if }k=i+1,\\f_k &\text{otherwise.}\end{cases}\]
Then $\sum \{f_i\} = \sum \{\widetilde f_i\} \mod I$ and the number of $i$ with $\widetilde f_i (v) > \widetilde f_{i+1}(v)$ for some $v\in X$ has decreased compared to the same quantity for the $f$'s, proving the claim.

From here, suppose 
\[
F=\sum_{i=1}^k\{f_i\}-\sum_{j=1}^\ell\{g_j\}
\] 
lies in the kernel of $\tau$.  Then, by Lemma~\ref{lemma:expfunctions}, there is an equality of multisets
\begin{equation}\label{eq:thing2}
\{f_i(v)\mid i=1,\dots,k\}=\{g_j(v)\mid j=1,\dots,\ell\}
\end{equation}
for each $v \in X$.  This implies that $k=\ell$, but it does not immediately imply that $f_i = g_i$ for each $i$; it could be the case, for instance, that $f_1(v) = g_2(v)$ and $f_2(v) = g_1(v)$.  However, the ``reordering of values'' claim above shows that, after adding an element of $I\subseteq\ker(\tau)$, we may assume that 
\begin{equation*}\label{eq:thing1}f_1(v) \le \cdots \le f_k(v)\text{ and } g_1(v) \le \cdots \le g_\ell(v) \text{ for all }v\in X.
\end{equation*}
Combining this with \eqref{eq:thing2} ensures that $f_i = g_i$ for all $i=1,\ldots,k$, so $F=0$, as desired.
\end{proof}

From here, the proof of Theorem~\ref{thm:K(X)presentation} is almost immediate.

\begin{proof}[Proof of Theorem~\ref{thm:K(X)presentation}]
By definition,
\[
\uPE(X)=\frac{\PE(X)}{\langle e^\ell-1\mid \ell\in \L(X)\rangle}.
\]
Since
\[
\langle e^\ell - 1 \mid \ell \in \L(X)\rangle = \Z\big\{ e^{f+\ell} - e^f \mid f \in \PL(X), \; \ell \in \L(X)\big\},
\]
the presentation for $\PE(X)$ given by Proposition~\ref{prop:K(X)presentation} implies that
\[
\uPE(X) \cong \frac{\Z\{\PL(X)\}}{\Z\bigg\{\substack{\vspace{5bp}\{f\} + \{g\} - \{\max(f,g)\} - \{\min(f,g)\},\\     \{f+\ell\} - \{f\}}\;\Big|\;\substack{\vspace{5bp}f,g\;\in\;\PL(X),\\     \ell\;\in\; \L(X)}\bigg\}}.\qedhere
\]
\end{proof}

We note that one can also write down an explicit presentation for $\uPE(\Sigma)$ for any unimodular fan $\Sigma$ as a quotient of the Stanley--Reisner ring of $\Sigma$ \cite{CCKR2}. Since this presentation is not necessary for the results of this paper, we do not include it here.

Expanding on the discussion of the introduction, we now close this section by describing two closely-related motivations for the study of integral piecewise-exponential functions.

\subsection{Motivation 1: Polytope algebras}
\label{subsec:polytopealg}

Polytope algebras have been studied extensively over the last several decades; see, for example, the work of McMullen \cite{McMullen}, Khovanski\u{\i}--Pukhlikov \cite{PukhlikovKhovanskii}, and Morelli \cite{Morelli} for foundational work from various different perspectives, and a much earlier paper of Groemer \cite{Groemer}, which implies that the various perspectives are equivalent. In this subsection, we briefly discuss the relationship between the ring of integral piecewise-exponential function on $N_\R$ and the polytope algebra of lattice polytopes in $M_\R$. 

Following McMullen \cite{McMullen}, the {\bf polytope algebra of $M_\R$}, denoted $\Pi(M_\R)$, is defined as the free abelian group generated by (convex) polytopes $P\subseteq M_\R$ subject to two types of additive relations:
\begin{equation}\label{eqn:pat}
[P+u]=[P]
\end{equation}
for all polytopes $P\subseteq M_\R$ and $u\in M_\R$, and
\begin{equation}\label{eqn:paie}
[P]+[Q]=[P\cup Q]+[P\cap Q]
\end{equation}
for all polytopes $P,Q\subseteq M_\R$ such that $P\cup Q$ is convex. The group $\Pi(M_\R)$ is then endowed with a ring structure given by Minkowski summation of polytopes:
\[
[P]\cdot[Q]=[P+Q].
\]
We are most interested in the subalgebra $\Pi(M)\subseteq\Pi(M_\R)$ generated by lattice polytopes---those polytopes whose vertices lie in $M$. Upon observing that the intersection of two lattice polytopes is a lattice polytope so long as their union is convex, it follows that the additive relations among the generators of $\Pi(M)$ are the same as \eqref{eqn:pat} and \eqref{eqn:paie} after imposing the extra restrictions that $P$ and $Q$ are lattice polytopes and $u\in M$.

Given any lattice polytope $P\subseteq M_\R$, consider its {\bf support function} $h_P$, defined by
\begin{align*}
h_P:N_\R&\rightarrow\R\\
v&\mapsto \max_{u\in P}\langle u,v\rangle.
\end{align*}
Because the vertices of $P$ are lattice points in $M$, it follows that the support function $h_P$ is integral piecewise-linear on the normal fan $\Sigma_P$.  Since $\Sigma_P$ is complete, this implies that $h_P\in\PL(N_\R)$. 
 Moreover, every \emph{convex} function $f\in\PL(N_\R)$ arises as the support function of some lattice polytope $P\subseteq N_\R$.

One can now define a function
\begin{align}
\label{eq:PitoPE}
    \Pi(M)&\rightarrow\uPE(N_\R)\\
\nonumber    [P]&\mapsto[e^{h_P}].
\end{align}
To see that this is well-defined under the defining relations \eqref{eqn:pat} and \eqref{eqn:paie} of $\Pi(M)$, note that
\begin{itemize}
\item $[e^{h_{P+u}}] = [e^{h_P}]$ (this follows from the fact that $h_{P+u} = h_P + u$ and $[e^u]=1 \in \uPE(N_\R)$); and
\item $[e^{h_P}] + [e^{h_Q}] = [e^{h_{P \cup Q}}] + [e^{h_{P \cap Q}}]$ when $P \cup Q$ is convex (this follows from Lemma~\ref{lemma:expfunctions} and the fact that $h_{P \cup Q} = \max(h_P, h_Q)$ whereas $h_{P \cap Q} = \min(h_P,h_Q)$).
\end{itemize}
Furthermore, the function \eqref{eq:PitoPE} is a ring homomorphism, as one checks by verifying that \[h_{P+Q} = h_P + h_Q.\]
In fact, this map is an isomorphism.  (While neither of the cited results below proves exactly this statement, it essentially follows from either.)

\begin{theorem}[\cite{Brion} Proposition~5.1, \cite{PukhlikovKhovanskii} Proposition 5]\label{thm:polytopealgebra}
The function \eqref{eq:PitoPE} is a ring isomomorphism $\Pi(M) \rightarrow \uPE(N_\R)$.
\end{theorem}

There is a natural $\Z$-algebra homomorphism on $\Pi(M)$, which we call the {\bf Ehrhart function}, defined by
\begin{align*}
    \Pi(M)&\rightarrow\Z\\
    [P]&\mapsto |P\cap M|.
\end{align*}
One perspective on the motivation for our work is to study generalizations of the Ehrhart function in the setting of incomplete rational fans.

\begin{motivatingquestion}\label{mq:ehrhart}
    For which rational fan supports $X\subseteq N_\R$ does there exist a ``natural'' $\Z$-algebra homomorphism
    \[
    \Chi\colon \uPE(X)\rightarrow\Z
    \]
    generalizing the Ehrhart function in the case where $X=N_\R$?
\end{motivatingquestion}

The requirement that $\Chi$ is ``natural'' can be viewed, for now, as simply the imprecise desideratum that there is a general procedure for associating $\Chi$ to $X$.  The specific requirements that $\Chi$ will be required to satisfy will be made precise in Section~\ref{sec:Ehrhart} below.

\subsection{Motivation 2: $\K$-rings}

If $\Sigma$ is a unimodular fan in $N_\R$, then $\Sigma$ corresponds to a smooth toric variety $Z_\Sigma$ containing the torus $\T=\mathrm{Spec}(\C[M])$ \cite{Fulton,CoxLittleSchenck}.  The {\bf Grothendieck $\K$-ring of $Z_\Sigma$}, denoted $K(Z_\Sigma)$, is defined as the free abelian group generated by vector bundles on $Z_\Sigma$, subject to the additive relations
\[
[\mathcal{F}]=[\mathcal{E}]+[\mathcal{G}]
\]
whenever there is a short exact sequence
\[
0\rightarrow\mathcal{E}\rightarrow\mathcal{F}\rightarrow\mathcal{G}\rightarrow 0.
\]
Tensor product of vector bundles then endows this group with the structure of a ring.

In the special setting of smooth toric varieties, the $\K$-ring is spanned by line bundles \cite[Theorem~1.1]{BrionVergne}, which are equivalent to piecewise-linear functions on $\Sigma$. For some intuition on this equivalence, let $\mathcal{L}$ be a $\T$-equivariant line bundle on $Z_\Sigma$. Then the restriction of $\mathcal{L}$ to any torus orbit of $Z_\Sigma$ is geometrically trivial, and therefore amounts simply to a character of $\T$, which is an element of $M$. Thus, a $\T$-equivariant line bundle gives rise to an integral linear function on each cone of $\Sigma$, and the compatibility required of these linear functions is that they agree on faces. In other words, $\T$-equivariant line bundles naturally correspond to piecewise-linear functions, and two equivariant line bundles are isomorphic nonequivariantly if and only if their corresponding piecewise-linear functions differ by a linear function. These considerations lead to the well-known presentation (\cite[Theorems~4.2.1 and 4.2.12]{CoxLittleSchenck}) of the Picard group of $Z_\Sigma$---the group of line bundles up to isomorphism---in terms of piecewise-linear functions:
\[
\mathrm{Pic}(Z_\Sigma)=\uPL(\Sigma).
\]

Enhancing the above discussion of Picard groups to the more general setting of $\K$-rings, we recall the following important result relating $\K$-rings of smooth toric varieties and piecewise-exponential functions.

\begin{theorem}[\cite{BrionVergne} Theorem~2.4, \cite{Merkurjev} Theorem~4.3]\label{thm:krings}
    For any unimodular fan $\Sigma$, 
    \[
    K(Z_\Sigma)=\uPE(\Sigma).
    \]
\end{theorem}

\begin{remark}
    Brion and Vergne proved that the $\T$-equivariant $\K$-ring is equal to $\PE(\Sigma)$, while Merkurjev proved that the map to the ordinary $\K$-ring is the appropriate quotient; see also the discussion in Section 6 of \cite{VezzosiVistoli}.
\end{remark}

\begin{remark}
    More generally, if $\Sigma$ is merely rational but not necessarily unimodular, Anderson and Payne \cite{AndersonPayne} showed that $\PE(\Sigma)$ is naturally isomorphic to the $\T$-equivariant operational K-ring of the possibly singular toric variety $Z_\Sigma$.
\end{remark}

If $\Sigma$ is complete, then the $\K$-ring $K(Z_\Sigma)$ admits a $\Z$-linear map to the integers, called the {\bf holomorphic Euler characteristic}, defined by
\begin{align*}
    \Chi:K(Z_\Sigma)&\rightarrow\Z\\
    [\mathcal{F}]&\mapsto \sum_{i\geq 0}(-1)^i\dim H^i(Z_\Sigma,\mathcal{F}).
\end{align*}
A second perspective on the motivation for our work is to study generalizations of the holomorphic Euler characteristic in the setting of incomplete unimodular fans.

\begin{motivatingquestion}\label{mq:euler}
    For which unimodular fans $\Sigma$ does there exist a ``natural'' linear map
    \[
    \Chi:\uPE(\Sigma)\rightarrow\Z
    \]
    generalizing the holomorphic Euler characteristic in the case where $\Sigma$ is complete?
\end{motivatingquestion}

It is a well-known result of toric geometry that the Ehrhart function and the holomorphic Euler characteristic agree when $\Sigma$ is a complete unimodular fan \cite[Proposition~9.4.3]{CoxLittleSchenck}, so Motivating Questions \ref{mq:ehrhart} and \ref{mq:euler} are really different perspectives on the same question.

Much of our intuition and terminology in what follows comes from Motivating Question~\ref{mq:ehrhart}.  On the other hand, one of our inspirations for studying these questions at all---as described in the introduction---comes from a recent paper of Larson, Li, Payne, and Proudfoot \cite{LLPP}, which developed an analogue of the holomorphic Euler characteristic on $\K$-rings associated to (incomplete) Bergman fans of matroids, thereby providing a partial answer to Motivating Question~\ref{mq:euler} that we seek to generalize.

\section{Ehrhart fans}
\label{sec:Ehrhart}

In this section, we introduce the notion of ``Ehrhart fans.''  Before giving the precise definition, we begin by briefly recalling the motivation discussed in the introduction.

\subsection{Motivation}
Given any complete unimodular fan and convex function $f \in \PL(\Sigma)$, consider the lattice polytope
\begin{equation}
\label{eq:Pf}
P_f\coloneq \{m \in M_\R \; | \; m(u_\rho) \leq f(u_\rho) \;\; \forall \rho \in \Sigma(1)\}.
\end{equation}
Ehrhart theory is concerned with counting the number of lattice points in such polytopes---that is, with calculating $|P_f \cap M|$---and we now describe a fundamental recursion satisfied by these lattice-point counts. Each ray $\rho \in \Sigma(1)$ defines a face $P^\rho_f$ of $P_f$, given by
\[
P_f^\rho \coloneq \{m \in P_f \; | \; m(u_\rho) = f(u_\rho)\}.
\]
On the other hand, if $\delta_\rho \in \PL(\Sigma)$ is the {\bf Courant function}, which is the piecewise-linear function taking value $1$ at $u_\rho$ and value $0$ at all other ray generators, then as long as the function $f- \delta_\rho$ is still convex, we can define $P_{f-\delta_\rho}$ exactly as we defined $P_f$. In this situation, the lattice-point counts of $P_f$, $P_{f - \delta_\rho}$, and $P_f^\rho$ satisfy the following relation
\begin{equation}
    \label{eq:PfcapN}
    |P_f \cap M| = |P_{f-\delta_\rho} \cap M| + |P_f^\rho \cap M|,
\end{equation}
which simply reflects the fact that $P_{f-\delta_\rho}$ is obtained from $P_f$ by ``sliding" the face $P_f^\delta$ inward by one lattice step, as illustrated in Figure~\ref{fig:intro}.

Notice that the polytope $P_{f + m}$ is a lattice translation of $P_f$ whenever $m\in M$, and this translation has the same number of lattice points as the original polytope.  Thus, there is a well-defined function
\[
[f] \mapsto |P_f \cap M|
\]
on those elements of $\uPL(\Sigma)$ represented by convex piecewise-linear functions $f$.  The definition of Ehrhart fans is created to generalize this setting, encoding those fans $\Sigma$ that admit an integer-valued function on $\uPL(\Sigma)$ satisfying an analogue of the relation \eqref{eq:PfcapN}.

\subsection{Definition of Ehrhart fans}

The definition of Ehrharticity is inductive on dimension, relating $\Sigma$ to its star fans along all rays $\rho$. In preparation for the definition of Ehrharticity, we first describe a natural map between piecewise-linear functions on fans and piecewise-linear functions on star fans. We remind the reader that our definitions and conventions regarding star fans were given in Subsection~\ref{subsection:conventions}.

If $g\in \PL(\Sigma)$ vanishes on a cone $\sigma\in\Sigma$, then on each cone $\tau\in\Sigma^{(\sigma)}$ in the neighborhood of $\sigma$, the restriction of $g$ to $\tau$ is constant on fibers of the quotient map $N\rightarrow N^\sigma$.  Thus, $g$ yields a well-defined piecewise-linear function on the star fan, which we denote $g^\sigma\in\PL(\Sigma^\sigma)$. Using this construction, we define a homomorphism
\begin{align*}
    \Pic(\Sigma) &\rightarrow \Pic(\Sigma^\sigma)\\
    [f] &\mapsto [f]^{\sigma} \coloneq [(f - m)^\sigma],
\end{align*}
where $m \in M$ is any linear function that agrees with $f$ on $\sigma$; it is straightforward to see that $[(f-m)^\sigma] \in \Pic(\Sigma^\sigma)$ does not depend on the choice of $m$.

We are now prepared to introduce Ehrhart fans.

\begin{definition}\label{def:Ehrhart}
Let $\Sigma$ be a unimodular fan.  We say that $\Sigma$ is {\bf Ehrhart} if
\begin{enumerate}
    \item the star fan $\Sigma^\rho$ is Ehrhart for all $\rho \in \Sigma(1)$, and
    \item there exists a well-defined function $\Chi_\Sigma\colon \Pic(\Sigma) \rightarrow \Z$ such that $\Chi_\Sigma(0) = 1$ and 
    \[
    \Chi_\Sigma([f]) = \Chi_\Sigma([f-\delta_\rho]) + \Chi_{\Sigma^\rho}([f]^\rho)
    \]
    for all $f \in \PL(\Sigma)$ and all $\rho \in \Sigma(1)$.
\end{enumerate}
If $\Sigma$ is Ehrhart, we call the function $\Chi_\Sigma$ the {\bf Ehrhart polynomial} of $\Sigma$.
\end{definition}

The next result justifies the use of the terminology ``the Ehrhart polynomial" in Definition~\ref{def:Ehrhart}, as opposed to ``an Ehrhart function.''

\begin{proposition}
If $\Sigma$ is an Ehrhart fan, then there is a unique function $\Chi_\Sigma\colon\uPL(\Sigma)\rightarrow\Z$ satisfying the conditions of Definition~\ref{def:Ehrhart}. Furthermore, the induced function 
    \begin{align*}
    \widehat\Chi_\Sigma\colon\Z^{\Sigma(1)}&\rightarrow\Z\\
    f&\mapsto\Chi_\Sigma([f])
    \end{align*}
agrees with a polynomial in $\Q[x_\rho \; | \; \rho\in\Sigma(1)]$ of degree $\dim(\Sigma)$.
\end{proposition}

\begin{proof}
    The proof is by induction on dimension. When $\Sigma$ is $0$-dimensional,  $\Pic(\Sigma) = 0$ and thus the only function satisfying the conditions of Definition~\ref{def:Ehrhart} is the constant function $\Chi_\Sigma=1$, which agrees with a polynomial function of degree zero, verifying the base case.  
    
    Now suppose that the proposition is true for all Ehrhart fans of dimension less than $d$, and let $\Sigma$ be an Ehrhart fan of dimension $d$. Let $\Chi_\Sigma$ be an Ehrhart function of $\Sigma$---which may not, a priori, be unique or given by a polynomial. 
    
    We first prove that $\Chi_\Sigma$ is unique. To prove this, it suffices to observe that repeated applications of the Ehrhart recursion 
    \[
    \Chi_\Sigma([f]) = \Chi_\Sigma([f-\delta_\rho]) + \Chi_{\Sigma^\rho}([f]^\rho)
    \]
    uniquely determine the value of $\Chi_\Sigma([f])$ for any $f\in\PL(\Sigma)=\Z^{\Sigma(1)}$ from the Ehrhart functions $\Chi_{\Sigma^\rho}$ and the initial value $\Chi_\Sigma(0)=1$. Since the Ehrhart functions $\Chi_{\Sigma^\rho}$ are all unique, by the induction hypothesis, it follows that $\Chi_\Sigma$ is unique.
    
    Now consider the induced function
    \begin{align*}
    \widehat\Chi_\Sigma\colon\Z^{\Sigma(1)}&\rightarrow\Z\\
    f&\mapsto\Chi_\Sigma([f])
    \end{align*}
    By the induction hypothesis, for each $\rho\in\Sigma(1)$, there is a polynomial $P_\rho\in\Q[x_\tau \; | \; \tau\in\Sigma(1)]$ of degree $\dim(\Sigma^\rho)$ such that
    \[
    \Chi_{\Sigma^\rho}([f]^\rho)=P_\rho(f)\;\;\;\text{ for every }\;\;\;f\in\Z^{\Sigma(1)},
    \]
    where $P_\rho(f)$ denotes the evaluation of $x_\tau$ at $f(u_\tau)$ for all $\tau\in\Sigma(1)$. Using that integer-valued polynomials in $\Q[x_\tau \; | \; \tau\in\Sigma(1)]$ have a $\Z$-basis given by products of univariate polynomials of the form
    \[
    {x_\tau\choose \alpha}=\frac{x_\tau(x_\tau-1)\dots (x_\tau-\alpha+1)}{\alpha!}
    \]
    with $\alpha\geq 0$ \cite[Proposition~XI.1.12]{CahenChabert}, it follows that there are integers $c_\alpha^\rho\in\Z$ such that
    \[
    P_\rho=\sum_{\alpha\in\N^{\Sigma(1)}\atop \alpha_{\rho}>0} c_\alpha^\rho{x_\rho-1\choose\alpha_\rho-1}\prod_{\tau\neq\rho}{x_\tau\choose \alpha_\tau}.
    \]
    We now claim that the integers $c_\alpha^\rho$ are independent of $\rho$.

    Let us compare coefficients $c_\alpha^\rho$ and $c_\alpha^\tau$ for distinct rays $\rho,\tau\in\Sigma(1)$. Since the $\alpha$ in question appears in a subscript of $c_\alpha^\rho$, then $\alpha_\rho>0$ by definition, and similarly, $\alpha_\tau>0$. Since $P_\rho(f)=\Chi_{\Sigma^\rho}([f]^\rho)$ only depends on the values $f(u_\tau)$ for which $\tau$ is a ray in the neighborhood $\Sigma^{(\rho)}$, and similarly for $P_\tau(f)$, it follows that $c_\alpha^\rho=c_\alpha^\tau=0$ if $\rho$ and $\tau$ do not lie on a common two-dimensional cone of $\Sigma$. Thus, suppose that there is a two-dimensional cone $\sigma\in\Sigma(2)$ such that $\sigma(1)=\{\rho,\tau\}$. Then the Ehrhart recursion applied to $\Sigma^\rho$ and $\Sigma^\tau$ implies that
    \[
    P_\rho(f)-P_\rho(f-\delta_\tau)=\Chi_{\Sigma^\sigma}([f]^\sigma)=P_\tau(f)-P_\tau(f-\delta_\rho).
    \]
    Using that
    \[
    {x_\tau\choose \alpha_\tau}-{x_\tau-1\choose\alpha_\tau}={x_\tau-1\choose\alpha_\tau-1},
    \]
    this implies that
    \[
    \sum_{\alpha\in\N^{\Sigma(1)}\atop \alpha_{\rho},\alpha_\tau>0} c_\alpha^\rho{x_\rho-1\choose\alpha_\rho-1}{x_\tau-1\choose\alpha_\tau-1}\prod_{\tau'\neq\rho,\tau}{x_{\tau'}\choose \alpha_{\tau'}}=\sum_{\alpha\in\N^{\Sigma(1)}\atop \alpha_\tau,\alpha_{\rho}>0} c_\alpha^\tau{x_\tau-1\choose\alpha_\tau-1}{x_\rho-1\choose\alpha_\rho-1}\prod_{\tau'\neq\rho,\tau}{x_{\tau'}\choose \alpha_{\tau'}}
    \]
    Thus, we see that $c_\alpha^\rho=c_\alpha^\tau$ for all $\alpha$ with $\alpha_\rho,\alpha_\tau>0$, as desired. 
    
    Now for each $\alpha\in\N^{\Sigma(1)}\setminus\{0\}$, define $c_\alpha\coloneq c_\alpha^\rho$ where $\rho\in\Sigma(1)$ is any ray for which $\alpha_\rho>0$. By our argument in the previous paragraph, the definition of $c_\alpha$ does not depend on the choice of $\rho$. Define the polynomial
    \[
    P=1+\sum_{\alpha\in\N^{\Sigma(1)}\setminus\{0\}}c_\alpha\prod_{\tau\in\Sigma(1)}{x_\tau\choose \alpha_\tau}\in\Q[x_\tau \; | \; \tau\in\Sigma(1)],
    \]
    which has degree
    \[
    1+\max_\rho\dim(\Sigma^\rho)=d. 
    \]
    For any $f\in\Z^{\Sigma(1)}$ and $\rho\in\Sigma(1)$, we compute
    \[
    P(f)-P(f-\delta_\rho)=P_\rho(f)=\Chi_{\Sigma^\rho}([f]^\rho)=\widehat\Chi_\Sigma(f)-\widehat\Chi_\Sigma(f-\delta_\rho),
    \]
    and it follows that $P-\widehat\Chi_\Sigma$ is a constant. Since both $P$ and $\widehat\Chi_\Sigma$ take value one at $0$, the two functions are thus equal. Therefore, $\widehat\Chi_\Sigma=P$, completing the proof of the induction step.
\end{proof}

\subsection{Examples and nonexamples of Ehrhart fans}

Since every zero-dimensional fan is Ehrhart, we now turn to an investigation of one-dimensional Ehrhart fans.  The next result completely characterizes when a one-dimensional unimodular fan is Ehrhart.

\begin{proposition}\label{prop:onedimensional}
A one-dimensional unimodular fan $\Sigma$ is Ehrhart if and only if 
\[
\sum_{\rho\in\Sigma(1)}u_\rho=0.
\]
Moreover, if $\Sigma$ is Ehrhart, then for every $f\in\PL(\Sigma)$, we have
\[
\Chi_\Sigma([f])=1+\sum_{\rho\in\Sigma(1)}f(u_\rho).
\]
\end{proposition}

\begin{proof}
Let $\Sigma$ be a one-dimensional fan.  Then condition (1) of Ehrharticity is satisfied, because the star fans $\Sigma^\rho$ are all zero-dimensional, and thus Ehrhart.  Thus, it remains to study whether the second condition of Ehrharticity holds for $\Sigma$. Define the function
\begin{align*}
G_\Sigma:\PL(\Sigma)&\rightarrow\Z\\
f&\mapsto 1+\sum_{\rho\in\Sigma(1)}f(u_\rho).
\end{align*}
Notice that $G_\Sigma(0)=1$, and since the Ehrhart polynomial of $\Sigma^\rho$ is the constant function $1$ for every $\rho\in\Sigma(1)$, it follows that $G_\Sigma$ satisfies the Ehrhart recursion:
\[
G_\Sigma(f)-G_\Sigma(f-\delta_\rho)=1=G_{\Sigma^\rho}([f]^\rho)\;\;\;\text{ for every }\;\;\;\rho\in\Sigma(1).
\]
Thus, $\Sigma$ is Ehrhart if and only if $G_\Sigma$ descends to a well-defined function on the quotient $\uPL(\Sigma)$. But $G_\Sigma$ is well-defined on $\Pic(\Sigma)$ if and only if 
\[
\sum_{\rho\in\Sigma(1)}f(u_\rho)=\sum_{\rho\in\Sigma(1)}(f+m)(u_\rho)
\]
for all $m\in M$, which is satisfied if and only if
\[
\sum_{\rho\in\Sigma(1)}u_\rho=0.\qedhere
\] 
\end{proof}

A pure, unimodular fan $\Sigma$ of dimension $d$ is said to be \textbf{balanced} (or \textbf{tropical with uniform weights}) if, for every $\tau\in\Sigma(d-1)$, we have
\[
\sum_{\sigma\in\Sigma(d)\atop \tau\subseteq\sigma}u_{\sigma(1)\setminus\tau(1)}\in\mathrm{span}(\tau).
\]
Equivalently, notice that $\Sigma$ is balanced if and only if the ray generators of the one-dimensional star fans $\Sigma^\tau$ add to zero for every $\tau\in\Sigma(d-1)$. Since star fans of Ehrhart fans are Ehrhart, Proposition~\ref{prop:onedimensional} implies the following result.

\begin{corollary}\label{corollary:balanced}
    Every Ehrhart fan is balanced.
\end{corollary}

Corollary~\ref{corollary:balanced} provides a wealth of unimodular fans that are not Ehrhart: any unimodular fan that is not balanced. Proposition~\ref{prop:onedimensional} shows that the converse of Corollary~\ref{corollary:balanced} is true in dimension one, and one might naturally wonder if the converse remains true in higher dimensions. It does not, and in order to give a counterexample, we now characterize two-dimensional Ehrhart fans.

\begin{proposition}\label{prop:dimensiontwo}
Let $\Sigma$ be a pure, unimodular fan of dimension two. Then $\Sigma$ is Ehrhart if and only if there exist integers $a_\rho\in\Z$ for every $\rho\in\Sigma(1)$ such that
\begin{equation}\label{eqn:dim2prop1}
\sum_{\sigma\in\Sigma(2)\atop \rho\subseteq\sigma}u_{\sigma(1)\setminus\rho}=a_\rho u_\rho
\end{equation}
and
\begin{equation}\label{eqn:dim2prop2}
\sum_\rho (2-a_\rho)u_\rho = 0.
\end{equation}
Moreover, if $\Sigma$ is Ehrhart, then for every $f\in\PL(\Sigma)$, we have
\[
\Chi_\Sigma([f])=1+\sum_{\rho\in\Sigma(1)}(1-a_\rho/2)f(u_\rho)+\sum_{\sigma\in\Sigma(2)\atop \sigma(1)=\{\rho_1,\rho_2\}}f(u_{\rho_1})f(u_{\rho_2})-\sum_{\rho\in\Sigma(1)}\frac{a_\rho}{2} f(u_\rho)^2.
\]
\end{proposition}

\begin{proof}
The first condition of Ehrharticity requires that all one-dimensional star fans $\Sigma^\rho$ are Ehrhart, but by Proposition~\ref{prop:onedimensional}, this is equivalent to the existence of integers $a_\rho\in\Z$ for every $\rho\in\Sigma(1)$ such that \eqref{eqn:dim2prop1} holds. Assuming that this is the case, define the function
\begin{align*}
G_\Sigma:\PL(\Sigma)&\rightarrow\Z\\
f&\mapsto 1+\sum_{\rho\in\Sigma(1)}\big(1-\frac{a_\rho}{2}\big)f(u_\rho)+\sum_{\sigma\in\Sigma(2)\atop \sigma(1)=\{\rho_1,\rho_2\}}f(u_{\rho_1})f(u_{\rho_2})-\sum_{\rho\in\Sigma(1)}\frac{a_\rho}{2} f(u_\rho)^2.
\end{align*}
Notice that $G_\Sigma(0)=1$, and for any $\rho\in\Sigma(1)$, a direct computation verifies that
\[
G_\Sigma(f)-G_\Sigma(f-\delta_\rho)=1+\sum_{\sigma\in\Sigma(2)\atop\rho\subseteq\sigma}f(u_{\sigma(1)\setminus\rho})-a_\rho f(u_\rho)=\Chi_{\Sigma^\rho}([(f-m)^\rho]),
\]
where $m$ is any linear function such that $m(u_\rho)=f(u_\rho)$. Thus, $G_\Sigma$ satisfies the necessary properties of an Ehrhart polynomial, so $\Sigma$ is Ehrhart if and only if $G_\Sigma$ descends to a well-defined function on $\uPL(\Sigma)$. A direct computation then shows that the quadratic part of $G_\Sigma$ descends if and only if \eqref{eqn:dim2prop1} holds, while the linear part of $G_\Sigma$ descends if and only if \eqref{eqn:dim2prop2} holds.
\end{proof}

We now provide a counterexample to the converse of Corollary~\ref{corollary:balanced}. The fan in the next example first appeared in work of Babaee and Huh \cite{BabaeeHuh}.

\begin{example}
Define 14 rays $\rho_1,\dots,\rho_{14}$ in $\R^4$ as the nonnegative spans of the columns in the following array:
\[
\begin{array}{cccccccccccccc}
u_1 & u_2 & u_3 & u_4 & u_5 & u_6 & u_7 & u_8 & u_9 & u_{10} & u_{11} & u_{12} & u_{13} & u_{14}\\
\hline
1 & 0 & 0 & 0 & 0 & 0 & 0 & 0 & 1 & 1 & 1 &-1 &-1 &-1\\
0 & 1 & 0 & 0 &-1 & 0 & 0 & 1 & 0 & 1 &-1 & 0 &-1 & 1\\
0 & 0 & 1 & 0 & 0 &-1 & 0 & 1 &-1 & 0 & 1 & 1 & 0 &-1\\
0 & 0 & 0 & 1 & 0 & 0 &-1 & 1 & 1 &-1 & 0 &-1 & 1 & 0\\
\end{array}
\]\\
Construct a fan $\Sigma$ by connecting the 14 rays according to the following graph:
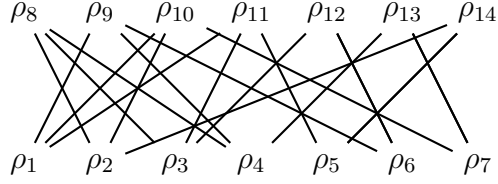
\begin{center}
\begin{tikzpicture}
\node (e1) at (1,0) {$\rho_1$};
\node (e2) at (2,0) {$\rho_2$};
\node (e3) at (3,0) {$\rho_3$};
\node (e4) at (4,0) {$\rho_4$};
\node (-e2) at (5,0) {$\rho_5$};
\node (-e3) at (6,0) {$\rho_6$};
\node (-e4) at (7,0) {$\rho_7$};
\node (f1) at (1,2) {$\rho_8$};
\node (f2) at (2,2) {$\rho_9$};
\node (f3) at (3,2) {$\rho_{10}$};
\node (f4) at (4,2) {$\rho_{11}$};
\node (-f2) at (5,2) {$\rho_{12}$};
\node (-f3) at (6,2) {$\rho_{13}$};
\node (-f4) at (7,2) {$\rho_{14}$};
\draw[thick] (e1) -- (f2);
\draw[thick] (e1) -- (f3);
\draw[thick] (e1) -- (f4);
\draw[thick] (e2) -- (f1);
\draw[thick] (e2) -- (f3);
\draw[thick] (e2) -- (-f4);
\draw[thick] (e3) -- (f1);
\draw[thick] (e3) -- (f4);
\draw[thick] (e3) -- (-f2);
\draw[thick] (e4) -- (f1);
\draw[thick] (e4) -- (f2);
\draw[thick] (e4) -- (-f3);
\draw[thick] (f2) -- (-e3);
\draw[thick] (f3) -- (-e4);
\draw[thick] (f4) -- (-e2);
\draw[thick] (-e2) -- (-f4);
\draw[thick] (-e3) -- (-f2);
\draw[thick] (-e4) -- (-f3);
\draw[thick] (-f2) -- (-e3);
\draw[thick] (-f3) -- (-e4);
\draw[thick] (-f4) -- (-e2);
\end{tikzpicture}
\end{center}
It follows from the work of Babaee and Huh that this is, in fact, a balanced unimodular fan. Moreover, we can compute the integers $a_i=a_{\rho_i}$ from this data, and we get the sequence
\[
\begin{array}{cccccccccccccc}
a_1 & a_2 & a_3 & a_4 & a_5 & a_6 & a_7 & a_8 & a_9 & a_{10} & a_{11} & a_{12} & a_{13} & a_{14}\\
\hline
3 & 3 & 3 & 3 & 0 & 0 & 0 & 1 & 1 & 1 & 1 & 0 & 0 & 0\\
\end{array}
\]\\
A direct computation now shows that
\[
\sum_{i=1}^{14}(2-a_i)u_i=(-4,-2,-2,-2)\neq 0,
\]
from which Proposition~\ref{prop:dimensiontwo} implies that $\Sigma$ is not Ehrhart.
\end{example}

\subsection{Leading terms of Ehrhart polynomials}

Given our motivation, one would hope that the Ehrhart polynomial of an Ehrhart fan reproduces some of the key properties of the classical Ehrhart polynomial $\Chi_P$ of a lattice polytope $P$ in $M_\R$, 
\begin{equation}
    \label{eq:chiP}
\Chi_P(t) = |tP \cap M|.
\end{equation}
One such property is that the leading coefficient of $\Chi_P$ is the volume of $P$. It is the goal of this subsection to prove a corresponding property for Ehrhart fans in general.  More specifically, for any Ehrhart fan $\Sigma$ of dimension $d$, we prove that the degree-$d$ piece of the Ehrhart polynomial $\Chi_\Sigma$ is equal to the ``volume polynomial'' of $\Sigma$, which has been studied in various other contexts, including \cite{Eur,NathansonRoss,RossLorentzian}. We begin by briefly recalling the definition of the volume polynomial.

For any unimodular fan $\Sigma$, in the ring $\text{Fun}(|\Sigma|, \R)$ of all functions $|\Sigma| \rightarrow \R$, let $\PP^*(\Sigma)$ denote the subring generated by Courant functions on $\Sigma$:
\[
\PP^*(\Sigma)\coloneq\Z[\delta_\rho\mid\rho\in\Sigma(1)] \subseteq \text{Fun}(|\Sigma|, \R).
\]
An element of $\PP^*(\Sigma)$ is the restriction of a polynomial function on each cone $\sigma\in\Sigma$, so the elements of $\PP^*(\Sigma)$ are ``piecewise-polynomial'' functions on $\Sigma$. If $\rho_1,\dots,\rho_k$ do not lie on a common cone, then $\delta_{\rho_1}\cdots\delta_{\rho_k}=0$, and moreover, this type of vanishing generates all relations in $\PP^*(\Sigma)$ (see \cite{Billera}). Thus, we have
\[
\PP^*(\Sigma)=\frac{\Z[x_\rho\mid\rho\in\Sigma(1)]}{\J_1},
\]
where
\[
\J_1 \coloneq \left\langle x_{\rho_1} \cdots x_{\rho_k} \; \left| \; \substack{\rho_1, \ldots, \rho_k \text{ do not }\\ \text{ span a cone of } \Sigma}\right.\right\rangle.
\]
The homogeneity of $\J_1$ implies that $\PP^*(\Sigma)$ is naturally graded by polynomial degree, and the degree-one piece is $\PP^1(\Sigma)=\PL(\Sigma)$. We define 
\[
\uPP^*(\Sigma)\coloneq\frac{\PP^*(\Sigma)}{\L(\Sigma)},
\]
so that
\[
\uPP^*(\Sigma)=\frac{\Z[x_\rho\mid\rho\in\Sigma(1)]}{\J_1+\J_2},
\]
where
\[
\J_2 \coloneq \Big\langle \sum_{\rho \in \Sigma(1)} m(u_\rho)x_\rho \; \Big| \; m \in M \Big\rangle.
\]
In toric geometry, we note that $\uPP^*(\Sigma)$ is naturally identified with the Chow ring of the smooth toric variety associated to $\Sigma$ \cite{Brion3}. An important property of $\uPP^*(\Sigma)$ is that it is spanned by square-free monomials \cite[Proposition 5.5]{AHK}. In particular, it follows that $\uPP^k(\Sigma)=0$ for all $k>\dim(\Sigma)$.

Balanced fans---and tropical fans, more generally---play an especially important role in this story due to the existence of nontrivial linear degree maps
\[
\deg:\uPP^d(\Sigma)\rightarrow\Z
\]
on the top-graded piece of $\uPP^*(\Sigma)$; in toric geometry, this is the degree map on the dimension-zero Chow group, and it features in tropical intersection theory.  In the specific case of balanced fans, it is a consequence of the balancing equation that any balanced fan admits a degree map that is uniquely determined by the property that
\[
\deg(x_{\rho_1}\cdots x_{\rho_d})=1
\]
whenever $\rho_1,\dots,\rho_d$ are the rays of a cone $\sigma\in\Sigma(d)$ \cite[Proposition 5.6]{AHK}. The volume polynomial of a balanced fan is defined via this degree map.

\begin{definition}
Let $\Sigma$ be a balanced fan of dimension $d$. The \textbf{volume polynomial} of $\Sigma$ is the function
\begin{align*}
    V_\Sigma:\uPL(\Sigma)&\rightarrow\Z\\
    [f]&\mapsto \deg_\Sigma([f]^d),
\end{align*}
where we identify $\uPL(\Sigma)=\uPP^1(\Sigma)$ and carry out the multiplication $[f]^d$ within $\uPP^*(\Sigma)$.
\end{definition}

The name ``volume polynomial'' comes from the special case where $\Sigma$ is complete, in which a convex $f \in \PL(\Sigma)$ determines a lattice polytope $P_f$ as in \eqref{eq:Pf}, and $V_\Sigma([f])$ is the (appropriately normalized) volume of $P_f$; see, for example,~\cite[Chapter 13]{CoxLittleSchenck}.  The main result of this subsection is the following.

\begin{proposition}
\label{prop:volume}
Let $\Sigma$ be an Ehrhart fan of dimension $d$. Then  the leading term of the Ehrhart polynomial is $\frac{1}{d!}V_\Sigma$.
\end{proposition}

\begin{proof}
The proof is by induction on $d$. If $d=0$, then $\Chi_\Sigma$ and $V_\Sigma$ are both the constant function $1$, proving the base case. Now suppose that $\Sigma$ has positive dimension $d$ and that the proposition holds for all fans of dimension less than $d$.

By precomposing the Ehrhart polynomial and the volume polynomial with the quotient map $\PL(\Sigma)\rightarrow \uPL(\Sigma)$, we obtain polynomials $\widehat\Chi_\Sigma$ and $\widehat V_\Sigma$ from $\PL(\Sigma)=\Z^{\Sigma(1)}$ to $\Z$. Since a homogeneous polynomial is uniquely determined from its partial derivatives, we consider the partial derivatives of $\frac{1}{d!}\widehat V_\Sigma$ and the partial derivatives of the degree-$d$ homogeneous component of the inhomogeneous polynomial $\widehat\Chi_\Sigma$.

Denoting elements of $\PL(\Sigma)$ by $f=\sum_{\rho}z_\rho\delta_\rho$, the chain rule gives
\[
\frac{\partial}{\partial z_\rho}\widehat V_\Sigma(f)=\frac{\partial}{\partial z_\rho}\deg_\Sigma([f]^d)=d\cdot \deg_\Sigma([f]^{d-1}[\delta_\rho]).
\]
A standard computation (see, for example, Equation (2.6) in \cite{RossLorentzian}) then gives
\[
\deg_\Sigma([f]^{d-1}[\delta_\rho])=\deg_{\Sigma^\rho}(([f]^\rho)^{d-1}),
\]
so that 
\begin{equation}\label{eq.partial1}
\frac{\partial}{\partial z_\rho}\frac{1}{d!}\widehat V_\Sigma(f)=\frac{1}{(d-1)!}\widehat V_{\Sigma^\rho}((f-m)^\rho),
\end{equation}
where $m$ is any linear function that agrees with $f$ at $u_\rho$. On the other hand, the derivative with respect to $z_\rho$ of the leading term of $\widehat \Chi_\Sigma(f)$ is equal to the leading term of
\[
\widehat \Chi_\Sigma(f)-\widehat \Chi_\Sigma(f-\delta_\rho),
\]
which, by the Ehrhart condition, is equal to the leading term of $\widehat\Chi_{\Sigma^\rho}((f-m)^\rho)$. Thus, letting $\widehat\Chi_\Sigma^d$ denote the leading term of $\widehat\Chi_\Sigma$, we have
\begin{equation}\label{eq.partial2}
\frac{\partial}{\partial z_\rho}\widehat\Chi_\Sigma^d(f)=\widehat\Chi_{\Sigma^\rho}^{d-1}((f-m)^\rho).
\end{equation}
Equations \eqref{eq.partial1} and \eqref{eq.partial2}, along with the induction hypothesis, prove that $\frac{1}{d!}\widehat V_\Sigma$ and $\widehat\Chi_\Sigma^d$ have the same partial derivatives, and since they are each homogeneous polynomials, Euler's homogeneous function theorem then implies that they are equal.
\end{proof}

\begin{remark}
Once we prove in the next section that complete fans are Ehrhart and $\Chi_\Sigma$ is a lattice-point count (Theorem~\ref{thm:complete}), Proposition~\ref{prop:volume} can be viewed as a natural generalization of the fact that the leading coefficient of the Ehrhart polynomial \eqref{eq:chiP} is the volume of $P$.  In particular, fixing a lattice polytope $P$ and letting $\Sigma$ denote any unimodular refinement of its normal fan, we have $P = P_f$ for some convex function $f \in \PL(\Sigma)$, and Theorem~\ref{thm:complete} shows that the classical Ehrhart polynomial $\Chi_P$ is related to $\Chi_\Sigma$ by
\[\Chi_P(t) = \Chi_\Sigma([tf]).\]
The leading coefficient of $\Chi_P$ is thus given by evaluating the leading term of $\Chi_\Sigma$ at $[f]$, and this equals the volume of $P$ by Proposition~\ref{prop:volume} and the preceding discussion.
\end{remark}

\section{Complete fans are Ehrhart}\label{sec:complete}

The aim of this section is to prove that complete fans are Ehrhart and that the Ehrhart polynomial of a complete fan computes lattice points of corresponding polytopes. We remind the reader that $n$ denotes the dimension of our ambient vector space $N_\R$, which is also the dimension of any complete fan in $N_\R$. The following is the main result of this section.

\begin{theorem}\label{thm:complete}
Let $\Sigma$ be a complete unimodular fan.  Then $\Sigma$ is Ehrhart. Furthermore, for every convex $f\in\PL(\Sigma)$, 
\[
\Chi_\Sigma([f])=|P_f\cap M|,
\]
where $P_f$ is the lattice polytope
\[
P_f \coloneq \{m \in M_\R \; | \; m(u_\rho) \leq f(u_\rho) \;\; \forall \rho \in \Sigma(1)\},
\]
and for every concave $f\in\PL(\Sigma)$,
\[
\Chi_\Sigma([f])=(-1)^n|\mathrm{Int}(P_{-f})\cap M|,
\]
where $\mathrm{Int}(P)$ denotes the topological interior of the polytope $P$.
\end{theorem}

\begin{remark}
Since complete unimodular fans correspond to smooth complete toric varieties, one could use known results in toric geometry to prove Theorem~\ref{thm:complete}, where the role of the Ehrhart polynomial would be played by the holomorphic Euler characteristic of line bundles. We take a different approach, modeled on techniques from Ehrhart theory; this more elementary proof does not require any tools from algebraic geometry. While the techniques below may be well-known to an expert in Ehrhart theory, it is useful to spell them out as we will generalize these ideas in the next section.
\end{remark}

Before embarking on the proof of Theorem~\ref{thm:complete}, we introduce notation for certain subfans of $\Sigma$ that will play a special role in the proof, and we state two key lemmas regarding these subfans (we hold off on the proofs of these lemmas until the end of this section).  In both lemmas, for $f \in\PL(\Sigma)$ and $m\in M$, let $\Sigma_{f,m}$ be the subfan of $\Sigma$ given by
\[
\Sigma_{f,m} \coloneq \{\sigma\in \Sigma \mid m(x) \le f(x)\text{ for all }x\in\sigma\}.
\]

\begin{lemma}\label{lem:finite}
    Let $\Sigma$ be a complete unimodular fan. For any $f\in\PL(\Sigma)$, there are only finitely many $m\in M$ such that
    \[
    \sum_{\sigma\in \Sigma_{f,m}}(-1)^{\dim \sigma}\neq 0.
    \]
\end{lemma}

\begin{lemma}\label{lem:convex}
    Let $\Sigma$ be a complete unimodular fan. For any convex $f\in\PL(\Sigma)$,
    \[
    \sum_{\sigma\in \Sigma_{f,m}}(-1)^{\dim \sigma}=\begin{cases}
        (-1)^n & m\in P_f\\
        0 & \text{else,}
    \end{cases}
    \]
    and for any concave $f\in\PL(\Sigma)$,
    \[
    \sum_{\sigma\in \Sigma_{f,m}}(-1)^{\dim \sigma}=\begin{cases}
        1& m\in \mathrm{Int}(P_{-f})\\
        0 & \text{else.}
    \end{cases}
    \]
\end{lemma}

We now prove Theorem~\ref{thm:complete} assuming the two lemmas above.

\begin{proof}[Proof of Theorem~\ref{thm:complete}]
Let $\Sigma$ be a complete unimodular fan of dimension $n$, and define a function $\Chi_\Sigma\col \PL(\Sigma)\to \Z$ by
\[
\Chi_\Sigma(f) = (-1)^n \sum_{m\in M} \left(\sum_{\sigma\in \Sigma_{f,m}}(-1)^{\dim \sigma}\right)\in \Z.
\]
Note that the outer sum is finite by Lemma~\ref{lem:finite}, so $\Chi_\Sigma(f)$ is a well-defined integer. Given any $m'\in M$, observe that
\begin{eqnarray*}
\Chi_\Sigma(f+m') &=& (-1)^n \sum_{m\in M} \sum_{\sigma\in \Sigma_{f\!+\!m',m}} (-1)^{\dim \sigma} \\ &=& (-1)^n \sum_{m\in M} \sum_{\sigma\in \Sigma_{f\!+\!m',m\!+\!m'}} (-1)^{\dim \sigma} \\ &=& \Chi_\Sigma(f),
\end{eqnarray*}
where the final equality uses that $\Sigma_{f,m} = \Sigma_{f+m', m+m'}$ for all $m'\in M$, which follows directly from the definition of $\Sigma_{f,m}$. Hence, $\Chi_\Sigma$ gives a well-defined function $\Chi_\Sigma:\Pic(\Sigma)\to\Z$. Moreover, for any convex $f\in\PL(\Sigma)$, Lemma~\ref{lem:convex} implies that
\begin{equation}\label{eq:countinglatticepoints}
\Chi_\Sigma([f])=(-1)^n \sum_{m\in P_f}(-1)^n=|P_f\cap M|,
\end{equation}
and for any concave $f\in\PL(\Sigma)$, Lemma~\ref{lem:convex} implies that
\[
\Chi_\Sigma([f])=(-1)^n \sum_{m\in \mathrm{Int}(P_{-f})}1=(-1)^n|\mathrm{Int}(P_{-f})\cap M|.
\]
It remains to prove that $\Sigma$ is Ehrhart with respect to the function $\Chi_\Sigma$.

We prove that $\Sigma$ is Ehrhart by induction on $n$.  First, if $\Sigma$ is complete then so is $\Sigma^\rho$ for all $\rho\in \Sigma(1)$.  Since $\dim \Sigma^\rho = \dim \Sigma - 1$, it follows inductively that $\Sigma^\rho$ is Ehrhart, verifying Property (1) of Ehrharticity.

In order to verify Property (2), we first note that
\[
\Chi_\Sigma(0)=|P_0\cap M|=1,
\]
where the first equality is a special case of \eqref{eq:countinglatticepoints} and the second equality is because $P_0=\{0\}.$ Lastly, we must show that, for any ray $\rho\in \Sigma(1)$, the function $\Chi_\Sigma$ satisfies the recursion
\[
\Chi_\Sigma([f]) - \Chi_\Sigma([f-\delta_\rho]) = \Chi_{\Sigma^\rho} ([f]^\rho).
\]
The left-hand side above is
\begin{equation}\label{eq:thingie}
    (-1)^n \sum_{m\in M}\left(\sum_{\sigma} (-1)^{\dim \sigma}\right),
\end{equation}
where $\sigma$ ranges over cones in $\Sigma_{f,m} \setminus \Sigma_{f-\delta_\rho,m}$.  Note that the expression in parentheses is zero outside a finite subset of $M$, since it is a difference of two expressions enjoying the same property.  Using the definitions of $\Sigma_{f,m}$ and $\Sigma_{f-\delta_\rho,m}$, the expression~\eqref{eq:thingie} is
\begin{equation}\label{eq:thingie2}
(-1)^n \sum_{\substack{m\in M \\ m(u_\rho) =f(u_\rho)}} \sum_{\substack{\sigma \in \Sigma, \rho \preceq \sigma, \\ m(x) \le f(x) \\ \text{for all } x\in \sigma}} (-1)^{\dim \sigma}.
\end{equation}
Now the cones in $\Sigma$ containing $\rho$ are exactly the cones of the star fan $\Sigma^\rho$, so this expression can be rewritten in terms of $\Sigma^\rho$. In order to do so, let $m_0\in M$ be any element such that $m_0(u_\rho) = f(u_\rho)$, and substitute $m' = m-m_0$.  Then~\eqref{eq:thingie2} is
\begin{eqnarray*}
(-1)^{n-1} \sum_{\substack{m'\in M \\  m'(u_\rho) =0 }} \, \sum_{\substack{\sigma \in \Sigma^\rho,  \\ m'(x) \le (f-m_0)(x) \\ \text{for all } x\in \sigma}} (-1)^{\dim \sigma}, 
\end{eqnarray*}
and this expression is exactly $\Chi_{\Sigma^\rho} ([f]^\rho)$. This completes the proof of Theorem~\ref{thm:complete}, modulo the proofs of Lemmas~\ref{lem:finite} and \ref{lem:convex}, which we provide below.
\end{proof}

In order to prove Lemmas~\ref{lem:finite} and \ref{lem:convex}, it will be helpful to reinterpret the sum 
\[
\sum_{\sigma\in\Sigma_{f,m}}(-1)^{\dim\sigma}
\]
in terms of compactly supported Euler characteristics. For a topological space $X$, let 
\[
\Chi_c(X) \coloneq \sum_i (-1)^i \dim H^i_c(X;\Q)\in\Z
\] 
denote the compactly supported Euler characteristic of $X$. The essential properties we require of $\Chi_c$ are the following:
\begin{enumerate}
    \item $\Chi_c$ is additive: for $Z\subseteq X$ a closed subspace of a sufficiently nice space $X$, then
    \[
    \Chi_c(X)=\Chi_c(Z)+\Chi_c(X\setminus Z).
    \]
    \item $\Chi_c(\RR^d) = (-1)^d$ for any $d\ge 0.$
\end{enumerate}
It follows that
\[
\Chi_c (\Sigma) = \sum_{\sigma\in\Sigma} \Chi_c(\mathrm{relint}(\sigma)) = \sum_{\sigma\in \Sigma}(-1)^{\dim \sigma}.
\]
Thus, we can replace the sums appearing in Lemmas~\ref{lem:finite} and \ref{lem:convex} with $\Chi_c(\Sigma_{f,m})$.\footnote{There is a slight abuse of notation; we should be writing $\Chi_c(|\Sigma|)$ rather than $\Chi_c(\Sigma)$. When making topological arguments below, we do not make a clear distinction between a polyhedral complex and its support, and it should be clear from context whether we are talking about the complex or the support.} With this, we are now prepared to prove the lemmas.

\begin{proof}[Proof of Lemma~\ref{lem:finite}]
Given any subfan $\Sigma'\subseteq\Sigma$, define
\[
V\coloneq\{m\in M_\R \mid \Sigma_{f,m} = \Sigma'\}\subseteq M_\R.
\]
We prove that $\Chi_c(\Sigma')=0$ whenever $V$ is unbounded. Since there are finitely-many subfans of $\Sigma$, this implies that there are finitely-many $m\in M$ such that $\Chi_c(\Sigma_{f,m})\neq 0$. Toward this end, suppose that $V$ is unbounded;  we shall show that the link of $\Sigma'$ is contractible, implying that $\Chi_c(\Sigma') = 0$.  We begin by studying the structure of the set $V$.  

Given $m\in M_\R$, the condition on $m$ that $\Sigma_{f,m}=\Sigma'$ is equivalent to the condition that for all ray generators $u_\rho$ of rays $\rho\in \Sigma(1)$, 
\[
\langle m, u_\rho\rangle \le f(u_\rho)\text{ if and only if }\rho \in \Sigma'.
\]
Therefore, $V$ is the intersection of the closed halfspaces
\[
\{m \in M_\R \mid \langle m,u_\rho\rangle \le f(u_\rho)\}
\]
for each ray generator $u_\rho$ of $\rho \in \Sigma'(1)$, and the open halfspaces
\[
\{m \in M_\R \mid \langle m,u_\rho\rangle > f(u_\rho)\}
\]
for each ray generator $u_\rho$ of $\rho\in\Sigma(1)\setminus\Sigma'(1)$. In particular, $V$ is an intersection of finitely many closed and open half-spaces.  Thus, if $V$ is unbounded, as is assumed, then it must be unbounded in a particular linear direction; that is, there exists $v\in M_\R\setminus \{0\}$ such that $v+ V \subseteq V$.

Now let $H$ be the hyperplane in $N_\R$ defined by
\[
H = \{ x \in N_\R\mid \langle v,x\rangle = 0\},
\] 
and let 
\[
H^+ = \{ x \in N_\R \mid \langle v,x\rangle >0\}\;\;\;\text{ and }\;\;\; H^- = \{ x \in N_\R \mid \langle v,x\rangle <0\}
\]
be the open half-spaces bounded by $H$.  We claim that whether a ray $\rho\in \Sigma(1)$ is in $\Sigma'(1)$ is almost entirely determined by its position relative to the hyperplane $H$.  The precise claim is as follows: let $u_\rho$ denote the primitive ray generator of $\rho \in \Sigma(1)$.  Then we claim that 
if $u_\rho \in H^-$ then $u_\rho\in \Sigma'$, and if $u_\rho\in H^+$ then $u_\rho\not \in \Sigma'$.  (If $u_\rho \in H$ then we cannot conclude either way, a priori.)  

To prove this claim, suppose first that $u_\rho\in H^-$, so $\langle v,u_\rho\rangle < 0$.  We claim that $\rho\in\Sigma'(1)$. If not, then $V\subseteq \{m\in M_\R: \langle m,u_\rho\rangle > f(u_\rho)\}$. Since $V$ is unbounded in the direction $v$, there exists $v_0\in V$ such that $\langle v_0+\lambda v, u_\rho \rangle > f(u_\rho)$ for all $\lambda\ge 0$.  But this contradicts that $\langle v, u_\rho\rangle < 0$.  A similar argument shows that if $u_\rho \in H^+$ then $\rho\not\in\Sigma'(1)$, proving the claim. The claim implies, in particular, that all rays of $\Sigma'$ are contained in $H\cup H^-$, and it follows that $\Sigma'\subseteq H\cup H^-$.

We now prove that the link of $\Sigma'$ is contractible. Here it will be convenient to use the ``combinatorial'' link of the unimodular fan $\Sigma'$, which is the simplicial complex $L\Sigma'$ glued from the simplices 
\[
L\sigma = \mathrm{conv}\{u_\rho~|~\rho\in\sigma(1)\}
\] 
along faces. Consider the space 
\begin{equation}\label{eq:Y-space}
Y \coloneq (L\Sigma' \cap H) \cup (L\Sigma\cap H^-).
\end{equation}
The space $L\Sigma\cap H^-$ is homeomorphic an open disk, and the whole space $Y$ is then homeomorphic to a partial closure of this open disk, which is contractible.  Note that $L\Sigma'$ is a subset of $Y$, since $\Sigma' \subseteq H \cup H^-$ as previously shown.  We shall construct a deformation retract from $Y$ to $L\Sigma'$.

We first set notation for certain deformation retracts of subsets of simplices. If $T$ is a finite set, write $\Delta^T$ for the standard simplex of dimension $|T|-1$ in $\R^T$, with vertices labeled by $T$.  If $S$ and $T$ are two nonempty sets then there is a canonical homeomorphism of the simplex $\Delta^{S\sqcup T}$ with the topological join $\Delta^S \ast \Delta^T$, defined as the quotient of $\Delta^S \times \Delta^T \times [0,1]$ by the equivalence relation $\sim$ generated by 
\[
(x,y,0) \sim (x',y,0) \text{ and } (x,y,1) \sim (x,y', 1)
\]
whenever $x,x'\in \Delta^S,$ and $y,y'\in \Delta^T$.  Thus there are canonical embeddings $\Delta^S \to \Delta^{S\sqcup T}$ and $\Delta^T \to \Delta^{S\sqcup T}$. We denote these images $\Delta^S$, respectively $\Delta^T$; they are closed faces of $\Delta^{S\sqcup T}$. 

The identification of $\Delta^{S\sqcup T}$ with the join $\Delta^S\ast\Delta^T$ makes evident a  deformation retraction of $\Delta^{S\sqcup T} \setminus \Delta^S$ onto $\Delta^T$.  In coordinates, it is the map
\begin{equation}\label{eq:in-coords}
\Delta^{S\sqcup T} \setminus \Delta^S\times[0,1]\to\Delta^{S\sqcup T} \setminus \Delta^S, \quad ((x,y,a),b)\mapsto (x,y,a(1-b)).
\end{equation}
It is important to note that these deformation retracts glue along faces: if $S'\subseteq S$ and $T'\subseteq T$ are nonempty subsets, then the following diagram commutes:
\begin{equation}\begin{aligned}\label{eq:commutes}
{\xymatrix@R=8mm@C=8mm{
\Delta^{S'\sqcup T'}\times [0,1]\ar[d]\ar[r] & \Delta^{S'\sqcup T'}\ar[d] \\ \Delta^{S\sqcup T}\times [0,1]\ar[r] &\Delta^{S\sqcup T}.
}}.
\end{aligned}\end{equation}

Now we return to the main proof. For each cone $\sigma\in \Sigma$ that meets $H^-$, let $\sigma'$ be the face of $\sigma$ spanned by its rays in $\Sigma'$, and let $\sigma''$ be the face of $\sigma$ spanned by its rays not in $\Sigma'$.  Then $L\sigma \cong L\sigma' \ast L\sigma''$.  Moreover, the subset
\[
Y_\sigma \coloneq (L\sigma'\cap H)\cup (L\sigma \cap H^-) \subseteq L\sigma \setminus L\sigma''
\]
admits a deformation retraction onto $L\sigma'$. Indeed, the deformation retraction of $L\sigma \setminus L\sigma''$ onto $L\sigma'$, given in~\eqref{eq:in-coords}, restricts to a deformation retraction of $Y_\sigma$ onto $L\sigma'$.  That is because, under the identification of $L\sigma$ with $L\sigma' \ast  L\sigma'' = (L\sigma' \times L\sigma'' \times [0,1])/\sim$, the set $Y_\sigma$ has the property that if $(x,y,a_1)\in Y_\sigma$ then $(x,y,a_2)\in Y_\sigma$ for all $a_2 \in [ a_1, 1]$.

These deformation retractions are compatible along faces by~\eqref{eq:commutes}, and so glue to give a deformation retraction of $Y$ onto $L\Sigma'$.  Therefore $L\Sigma'$ is contractible since $Y$ is contractible, and we conclude that $\Chi_c(L\Sigma') = 0$, as desired, concluding the proof of Lemma~\ref{lem:finite}.
\end{proof}

\begin{proof}[Proof of Lemma~\ref{lem:convex}]
We focus on the statement regarding convex functions; using that the negation of any concave function is convex, the statement regarding concave functions can then be obtained via appropriate modifications of the arguments. Suppose that $f\in\PL(\Sigma)$ is convex, so that
\[
f(x)+f(y)\geq f(x+y)\;\;\;\text{ for all }\;\;\;x,y\in N_\R.
\]
If $m\in P_f$, then it follows from definition that $\Sigma_{f,m}=\Sigma$, and we compute
\[
\Chi_c(\Sigma_{f,m})=\Chi_c(\R^n)=(-1)^n.
\]
This proves one case of the statement in the lemma; it remains to prove that $\Chi_c(\Sigma_{f,m})=0$ when $m\notin P_f$.

Assume that $m\in M\setminus P_f$, and define the set 
\[
U \coloneq \{x\in N_\R \mid m(x) > f(x)\}.
\]
The convexity of the piecewise-linear function $f$ implies that $U$ is an open convex polyhedral cone in $N_\R$, and the assumption that $m\notin P_f$ further implies that $U$ is nonempty. In particular, defining $Y$ as the intersection of $L\Sigma$ with the complement of $U$
\[
Y\coloneq L\Sigma\cap U^c,
\]
it follows that $Y$ is contractible.

By definition, a cone $\sigma\in\Sigma$ lies in $\Sigma_{f,m}$ if and only if its rays lie in $U^c$. In particular, $\Sigma_{f,m}\subseteq U^c$, implying that $L\Sigma_{f,m}\subseteq Y$. Moreover, an argument similar to that in the proof of Lemma~\ref{lem:finite} shows that $Y$ admits a deformation retract onto $L\Sigma_{f,m}$. The contractibility of $Y$ then implies the contractibility of $L\Sigma_{f,m}$, from which we conclude that $\Chi_c(\Sigma_{f,m})=0$, as desired.
\end{proof}

\section{Independence of fan structure}
\label{sec:independence}

In this section, we prove that the property of being Ehrhart is independent of the fan structure of $\Sigma$, instead depending only on the support of $\Sigma$. An important consequence of this independence---discussed at the end of this section---is that Ehrhart polynomials induce well-defined linear functions on rings of piecewise-exponential functions. 

We prove our main independence result via a theorem of W{\l}odarczyk \cite{wlodarczyk}---stated precisely as Theorem~\ref{thm:wlodarczyk} below---which allows us to connect any two unimodular fans by a sequence of stellar subdivisions and their inverses. In order to employ W{\l}odarczyk's result, we require two preliminary steps: the first describes how Ehrharticity interacts with fiber bundles of fans, and the second describes how Ehrharticity interacts with stellar subdivisions.

\subsection{Ehrharticity and fiber bundles}

In order to understand how Ehrharticity interacts with fiber bundles, we begin by recalling the definition of fiber bundles in the context of unimodular fans.

\begin{definition}\label{def:fiberbundle}
Let $0\to N' \to N \xrightarrow{\pi} N'' \to 0$ be a short exact sequence of finitely-generated free abelian groups, and let $\Sigma'$ and $\Sigma''$ be unimodular fans in $N'_\R$ and $N''_\R$, respectively.  (We view $N'_\R \subseteq N_\R$ by tensoring the short exact sequence with $\R$.)  A \textbf{$\Sigma'$-bundle over $\Sigma''$} is a unimodular fan $\Sigma$ in $N_\R$ such that
\begin{enumerate}
    \item $|\Sigma|\cap N_\R'=|\Sigma'|$;
    \item there exists a subfan $\widetilde\Sigma\subseteq\Sigma$ such that $\pi$ is injective on $\widetilde\Sigma$ and the function
    \begin{align*}
        \widetilde\Sigma&\rightarrow\Sigma''\\
        \widetilde\sigma&\mapsto\pi(\widetilde\sigma)
    \end{align*}
    is a bijection such that $\pi(\widetilde \sigma \cap N) = \pi(\widetilde \sigma) \cap N''$ for each $\widetilde \sigma \in \widetilde \Sigma$;
    \item the function
    \begin{align*}
        \Sigma'\times\widetilde\Sigma&\rightarrow\Sigma\\
        (\sigma',\widetilde\sigma)&\mapsto\sigma'+\widetilde\sigma
    \end{align*}
    is a bijection.
\end{enumerate}
\end{definition}

Combinatorially, the cones in a $\Sigma'$-bundle over $\Sigma''$ are the same as the cones of the product $\Sigma'\times\Sigma''$, but the linear structure of a bundle may be different from the linear structure of the product. For example, letting $\Sigma$ be the complete pointed one-dimensional fan, the image below depicts a $\Sigma$-bundle over itself which is not simply the product $\Sigma\times\Sigma$.

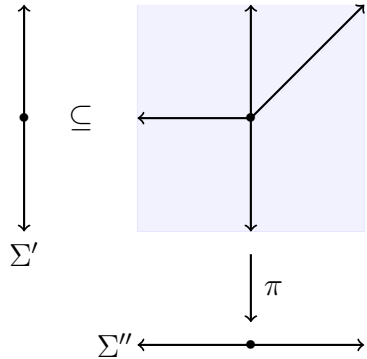
\begin{center}
\begin{tikzpicture}[scale=1.5]
\draw[draw=blue!10,fill=blue!10, fill opacity=.5]    (1,1) -- (1,-1) -- (-1,-1) -- (-1,1); 
\draw[thick,->] (0,0) -- (1,1); 
\draw[thick,->] (0,0) -- (0,1); 
\draw[thick,->] (0,0) -- (-1,0); 
\draw[thick,->] (0,0) -- (0,-1); 

\draw[thick,->] (-2,0) -- (-2,1); 
\draw[thick,->] (-2,0) -- (-2,-1); 

\draw[thick,->] (0,-2) -- (-1,-2); 
\draw[thick,->] (0,-2) -- (1,-2); 

\draw[] node at (0,0) {\scalebox{.7}{$\bullet$}};
\draw[] node at (-2,0) {\scalebox{.7}{$\bullet$}};
\draw[] node at (0,-2) {\scalebox{.7}{$\bullet$}};

\draw[] node at (-1.5,0) {$\subseteq$};
\draw[thick,->] (0,-1.2) -- (0,-1.8);
\draw[] node at (.2,-1.5) {$\pi$};

\draw[] node at (-2,-1.2) {$\Sigma'$};
\draw[] node at (-1.2,-2) {$\Sigma''$};
\end{tikzpicture}
\end{center}

\begin{remark}
    In the language of toric geometry, a $\Sigma'$-bundle over $\Sigma''$ corresponds to a smooth toric variety $Z_\Sigma$ along with a projection map to $Z_{\Sigma''}$ whose fibers are all isomorphic to $Z_{\Sigma'}$. This perspective informs but is not essential to what follows.
\end{remark}

The next result describes how Ehrharticity interacts with bundles having complete fibers.

\begin{proposition}\label{prop:fiberbundles}
Let $\Sigma'$ and $\Sigma''$ be unimodular fans. If $\Sigma'$ is complete, then a $\Sigma'$-bundle over $\Sigma''$ is Ehrhart if and only if $\Sigma''$ is Ehrhart.
\end{proposition}

\begin{proof}
Suppose that $\Sigma'$ is complete and that $\Sigma$ is a $\Sigma'$-bundle over $\Sigma''$. Throughout this proof, we use the notation of Definition~\ref{def:fiberbundle}.

First, let us assume that $\Sigma$ is Ehrhart. Choose any maximal cone $\sigma\in\Sigma'\subseteq\Sigma$. It follows directly from the definitions of fiber bundles and star fans that $\Sigma^\sigma\cong \Sigma''$. Since the star fan of $\Sigma$ at an arbitrary cone $\sigma$ can be obtained iteratively from $\Sigma$ by taking star fans at a ray, the assumption that $\Sigma$ is Ehrhart implies that $\Sigma''$ is Ehrhart.

Conversely, suppose that $\Sigma''$ is Ehrhart. We use induction on $\dim\Sigma$ to prove that $\Sigma$ is Ehrhart. In other words, we assume that a fiber bundle is Ehrhart whenever it has complete fibers, an Ehrhart base, and dimension less than $\dim\Sigma$. To check Condition (1) of Ehrharticity, let $\rho\in\Sigma(1)$. Then either $\rho\in\Sigma'(1)$ or $\rho\in\widetilde\Sigma(1)$. In the first case, it follows from the definitions that $\Sigma^\rho$ is a $(\Sigma')^\rho$-bundle over $\Sigma''$, and in the second case, $\Sigma^\rho$ is a $\Sigma'$-bundle over $(\Sigma'')^{\pi(\rho)}$. In either case, $\Sigma^\rho$ is Ehrhart by induction on $\dim\Sigma$.

To prove Condition (2) of Ehrharticity for the fiber bundle $\Sigma$, we define an explicit function $\Chi_\Sigma:\PL(\Sigma)\rightarrow\Z$; to do so, we must first introduce additional notation. Given $f\in \PL(\Sigma)$, define $f'\in\PL(\Sigma')$ as the restriction of $f$ to $\Sigma'$, and define $f''\in\PL(\Sigma'')$ as the unique piecewise-linear function such that
\[
f''(u_{\pi(\rho)})=f(u_\rho)\;\;\;\text{ for every }\;\;\;\rho\in\widetilde\Sigma(1).
\]
Choose a splitting of the quotient map $\pi:N_\R\rightarrow N_\R''$, allowing us to write $N_\R=N_\R'\oplus N_\R''$. Given any $m'\in M'$---in other words, an integral linear map on $N_\R'$---we can use the above splitting to view $m'$ to an integral linear map on $N_\R$. For any $f\in\PL(\Sigma)$, define (using the notation from Section~\ref{sec:complete})
\[
\Chi_\Sigma(f)=(-1)^{\dim\Sigma'}\sum_{m'\in M'}\Big(\Big (\sum_{\sigma\in\Sigma'_{f',m'}}(-1)^{\dim(\sigma)}\Big)\Chi_{\Sigma''}([(f-m')''])\Big).
\]

We first prove that $\Chi_\Sigma$ descends to $\uPL(\Sigma)$. Given any $\ell\in M$, we compute
\begin{align*}
    \Chi_\Sigma(f+\ell)&=(-1)^{\dim\Sigma'}\sum_{m'\in M'}\Big(\sum_{\sigma\in\Sigma'_{f'+\ell',m'}}(-1)^{\dim(\sigma)}\Big)\Chi_{\Sigma''}([(f+\ell-m')''])\\
    &=(-1)^{\dim\Sigma'}\sum_{m'+\ell'\in M'}\Big(\sum_{\sigma\in\Sigma'_{f'+\ell',m'+\ell'}}(-1)^{\dim(\sigma)}\Big)\Chi_{\Sigma''}([(f+\ell-m'-\ell')''])\\
    &=\Chi_\Sigma(f),
\end{align*}
where the final equality uses the fact that $\Sigma'_{f'+\ell',m'+\ell'}=\Sigma'_{f',m'}$, along with the observation that $(\ell-\ell')''$ is linear on $\Sigma''$. Thus, $\Chi_\Sigma$ descends to a function $\Chi_\Sigma:\uPL(\Sigma)\rightarrow\Z$.

Next, we observe that 
\[
\Chi_\Sigma(0)=(-1)^{\dim\Sigma'}\sum_{m'\in M'}\Big(\sum_{\sigma\in\Sigma'_{0,m'}}(-1)^{\dim(\sigma)}\Big)\Chi_{\Sigma''}([(-m')''])=\Chi_{\Sigma''}(0)=1,
\]
where the second equality uses Lemma~\ref{lem:convex} to conclude that $m'=0$ is the only element of $M'$ contributing to the sum. To complete the proof that $\Sigma$ is Ehrhart, it now remains to verify that
\begin{equation}\label{eq:fiberbundlerecursion}
\Chi_\Sigma([f]) - \Chi_\Sigma([f-\delta_\rho]) = \Chi_{\Sigma^\rho}([f]^\rho)
\end{equation}
for all $f\in\PL(\Sigma)$ and $\rho\in\Sigma(1)$. There are two cases to consider, depending on whether $\rho\in\Sigma'(1)$ or $\rho\in\widetilde\Sigma(1)$.

If $\rho\in\Sigma'(1)$, then the left-hand side of \eqref{eq:fiberbundlerecursion} can be written as 
\[
(-1)^{\dim\Sigma'}\sum_{m'\in M'}\Big(\sum_{\sigma\in\Sigma'_{f',m'}\setminus\Sigma'_{f'-\delta_\rho,m'}}(-1)^{\dim(\sigma)}\Big)\Chi_{\Sigma''}([(f-m')'']).
\]
As in the proof of Theorem~\ref{thm:complete}, the cones $\sigma$ appearing in the sum above can be identified with the cones of $(\Sigma')^\rho$, and then interpreting $\Sigma^\rho$ as a $(\Sigma')^\rho$-bundle over $\Sigma''$, it then follows that this difference is equal to $\Chi_{\Sigma^\rho}([f]^\rho)$, verifying \eqref{eq:fiberbundlerecursion} in this case.

If $\rho\in\widetilde\Sigma(1)$, then the left-hand side of \eqref{eq:fiberbundlerecursion} can be written as 
\[
(-1)^{\dim\Sigma'}\sum_{m'\in M'}\Big(\sum_{\sigma\in\Sigma'_{f',m'}}(-1)^{\dim(\sigma)}\Big)\Big(\Chi_{\Sigma''}([(f-m')''])-\Chi_{\Sigma''}([(f-m')''-\delta_{\pi(\rho)}])\Big).
\]
Since $\Sigma''$ is Ehrhart, the final difference in the above expression is equal to 
\[
\Chi_{(\Sigma'')^{\pi(\rho)}}(([f-m')]^{\pi(\rho)}),
\]
and then the interpretation of $\Sigma^\rho$ as a $\Sigma'$-bundle over $(\Sigma'')^{\pi(\rho)}$ implies the validity of \eqref{eq:fiberbundlerecursion} in this case, finishing the proof.
\end{proof}

The key insight in the proof of Proposition~\ref{prop:fiberbundles} is that the perspective from Section~\ref{sec:complete} on the Ehrhart polynomial of a complete fan can be ``fibered'' over any other Ehrhart polynomial, giving a candidate for the Ehrhart polynomial of the fiber product that we have proven is correct.  On the other hand, when the fiber bundle is a \emph{product}, there is also a natural candidate for the Ehrhart polynomial---namely, the product of the Ehrhart polynomials of the two factors---even without assuming that either factor is complete.  This candidate is also correct, and since this statement will be needed in Section~\ref{sec:matroidfans} of the paper, we prove it here before proceeding.

\begin{proposition}
\label{prop:products}
Let $\Sigma_1$ and $\Sigma_2$ be Ehrhart fans.  Then $\Sigma_1 \times \Sigma_2$ is Ehrhart, with Ehrhart polynomial given by the product of the Ehrhart polynomials of the factors:
\[
 \Chi_{\Sigma_1\times\Sigma_2}([f_1 \times f_2])=\Chi_{\Sigma_1}([f_1])\times\Chi_{\Sigma_2}([f_2]).
 \]
\end{proposition}
\begin{proof}
Let $d_i=\dim(\Sigma_i)$ and  $\Sigma = \Sigma_1 \times \Sigma_2$.  We proceed by induction on $d_1 + d_2$.

The base case is when $d_1 +d_2 = 0$, and hence both fans as well as their product are zero-dimensional and therefore Ehrhart. For the induction step, suppose products of Ehrhart fans are Ehrhart with the desired Ehrhart polynomial whenever the sum of their dimensions is less than $d$, and suppose that $d_1 + d_2 = d$.  Given a ray $\rho\in\Sigma(1)$, either $\rho\in\Sigma_1(1)$ or $\rho\in\Sigma_2(1)$. In the first case, we have
\[
\Sigma^\rho=\Sigma_1^{\rho}\times\Sigma_2,
\]
and similarly for the second case. It then follows from the induction hypothesis and the assumption that both $\Sigma_1$ and $\Sigma_2$ are Ehrhart that $\Sigma^\rho$ is Ehrhart for every $\rho\in\Sigma(1)$. Given a piecewise-linear function $f= f_1 \times f_2$ on $\Sigma$, we define
\[
\Chi_\Sigma([f])=\Chi_{\Sigma_1}([f_1])\Chi_{\Sigma_2}([f_2]).
\]
To show that this function satisfies the Ehrhart relation, let $\rho \in \Sigma(1)$, and suppose without loss of generality that $\rho\in\Sigma_1(1)$ so that $f-\delta_\rho=(f_1-\delta_\rho) \times f_2$. We then compute
\begin{align*}
\Chi_\Sigma([f])-\Chi_\Sigma([f-\delta_\rho])&=\Chi_{\Sigma_1}([f_1])\Chi_{\Sigma_2}([f_2])-\Chi_{\Sigma_1}([f_1-\delta_\rho])\Chi_{\Sigma_2}([f_2])\\
&=\Chi_{(\Sigma_1)^\rho}([f_1]^\rho)\Chi_{\Sigma_2}([f_2])\\
&=\Chi_{\Sigma^\rho}([f]^\rho),
\end{align*}
where we use the assumption that $\Sigma_1$ is Ehrhart in the second equality and the induction hypothesis in the third.
\end{proof}

\subsection{Ehrharticity and stellar subdivisions}
Next, we require an understanding of how Ehrharticity interacts with stellar subdivisions, which will build on our discussion of fiber bundles in the previous section. Let us briefly recall the construction of stellar subdivisions and the key properties of them that we require.

\begin{definition}
Let $\Sigma$ be a unimodular fan and let $\tau\in\Sigma$. Define
\[
u_\tau=\sum_{\rho\in\tau(1)}u_\rho\;\;\;\text{ and }\;\;\;\rho_\tau=\R_{\geq 0}u_\tau.
\]
The \textbf{stellar subdivision of $\Sigma$ at $\tau$}, denoted $\Sigma_\tau$, is the fan 
\[
\Sigma_\tau \coloneq \{\sigma \in \Sigma \; | \; u_\tau \notin \sigma\} \;\; \cup \;\;\bigcup_{\substack{ u_\tau \in \sigma }} \;\; \bigcup_{\substack{\pi \preceq \sigma,\\u_\tau \notin \pi}} \pi + \rho_\tau,
\]
where in the two unions, both $\sigma$ and $\pi$ range over cones of $\Sigma$.
\end{definition}

In other words, to construct $\Sigma_\tau$ from $\Sigma$, we replace each cone $\sigma$ of $\Sigma$ containing $\tau$ by the collection of all cones spanned by $u_\tau$ and a disjoint face of $\sigma$. Note that the stellar subdivision of a unimodular fan along any cone is, itself, a unimodular fan.

\begin{remark}
    In the language of toric geometry, a stellar subdivision, as defined above, corresponds to a blowup of a smooth toric variety along a toric subvariety. This perspective informs but is not essential to what follows.
\end{remark}

In order to study Ehrharticity of stellar subdivisions, it will be important to understand star fans at rays of a stellar subdivision. The next lemma describes star fans of stellar subdivisions explicitly. Recall that $\Sigma^{(\tau)}$ denotes the neighborhood of $\tau$ in $\Sigma$ (see Subsection~\ref{subsection:conventions}).

\begin{lemma}\label{lem:starfans}
Let $\Sigma$ be a unimodular fan in $N_\R$ and $\tau\in\Sigma$. Given any ray $\rho$ in the stellar subdivision $\Sigma_\tau$, the star fan $(\Sigma_\tau)^\rho$ is described by one of the following three cases.
\begin{enumerate}
\item If $\rho\notin\Sigma^{(\tau)}(1)$, then
\[
(\Sigma_\tau)^\rho=\Sigma^\rho.
\]
\item If $\rho\in\Sigma^{(\tau)}(1)\setminus\{\rho_\tau\}$, then $\tau$ projects onto a cone $\overline\tau\in\Sigma^\rho$, and the star fan of $\Sigma_\tau$ at $\rho$ is the stellar subdivision of the star fan $\Sigma^\rho$ at $\overline\tau$:
\[
(\Sigma_\tau)^\rho=(\Sigma^\rho)_{\overline\tau}.
\]
\item If $\rho=\rho_\tau$, then $(\Sigma_\tau)^\rho$ is a $\Pi_{\dim\tau-1}$-bundle over $\Sigma^\tau$, where $\Pi_k$ denotes the normal fan of the standard $k$-simplex in $\R^{k+1}/\R$.
\end{enumerate}
\end{lemma}

\begin{proof}
To justify the description of the star fan in Case (1), note that, in this case, the neighborhood $\Sigma^{(\rho)}$---and thus the star fan $\Sigma^\rho$---is unaffected by  stellar subdivision at $\tau$. 

To justify Case (2), note that any cone $\sigma$ in the neighborhood $(\Sigma_\tau)^{(\rho)}$ is either a cone of $\Sigma^{(\rho)}$ (if $u_\tau\notin\sigma$), or it is a cone of the form
\[
\pi + \rho_\tau
\]
for some $\pi\in\Sigma^{(\rho)}$ such that $\pi$ and $u_\tau$ lie in a common cone of $\Sigma^{(\rho)}$ and $u_\tau\notin\pi$. Upon projecting these two types of cones to $N_\R^\rho$, we obtain the two types of cones comprising the stellar subdivision of $\Sigma^\rho$ at $\overline\tau$.

Lastly, we consider Case (3). In order to describe $(\Sigma_\tau)^{\rho_\tau}\subseteq N_\R^{\rho_\tau}$ as a fiber bundle, we describe the fiber and the base explicitly. Define $\Sigma'$ to be the image of $\Sigma_\tau\cap\tau$ in the quotient space $N_\R^{\rho_\tau}$, and define $N_\R'=\Span_\R(\Sigma')\subseteq N_\R^{\rho_\tau}$. It follows from the definition of stellar subdivisions that $\Sigma'=(\Sigma_\tau)^{\rho_\tau}\cap N_\R'$ is the normal fan of the standard simplex. Let $\widetilde\Sigma$ be the maximal subfan of $(\Sigma_\tau)^{\rho_\tau}$ that omits the rays of $\Sigma'$, or in other words, omits the images of rays of $\tau$. It then follows from the definition of star fans that the natural quotient map 
\[
\pi:N_\R^{\rho_\tau}\rightarrow N_\R^{\rho_\tau}/N_\R' = N_\R^\tau
\]
is an injection of the cones of $\widetilde\Sigma$ onto the cones of the star fan $\Sigma^\tau$. Moreover, by the definition of stellar subdivisions, the cones of $(\Sigma_\tau)^{\rho_\tau}$ are in natural bijection with pairs $(\sigma',\widetilde\sigma)\in\Sigma'\times\widetilde\Sigma$. Putting this all together, it follows that $(\Sigma_\tau)^{\rho_\tau}$ is a fiber bundle over $\Sigma''=\Sigma^\tau$, with fibers $\Sigma'=\Pi_{\dim\tau-1}$, completing the proof.
\end{proof}

The other important ingredient in our study of Ehrharticity of stellar subdivisions is the observation that, if $\Sigma$ is a unimodular fan and $\tau \in \Sigma$, then for a ray $\rho \in \Sigma$, there is an associated Courant function on $\Sigma$ as well as on $\Sigma_\tau$, and these are generally different functions on $|\Sigma| = |\Sigma_\tau|$.  We denote the Courant functions on these two fans by
\[\delta_\rho \in \PL(\Sigma) \;\;\; \text{ and } \;\;\; \tilde\delta_\rho \in \PL(\Sigma_\tau)\]
to distinguish them (where the latter is also defined when $\rho = \rho_\tau$); the relationship between them is that
\[\delta_\rho = \begin{cases} \tilde\delta_\rho & \text{ if } \rho \notin \tau(1),\\ \tilde\delta_\rho + \tilde\delta_{\rho_\tau} & \text{ if } \rho \in \tau(1),\end{cases}\]
as one verifies by evaluating both sides at all the ray generators of $\Sigma_\tau$.

We are now prepared to discuss the main technical result of this section, which shows that Ehrharticity is insensitive to taking stellar subdivisions.

\begin{proposition}\label{prop:blowups}
If $\Sigma$ is a unimodular fan and $\tau\in\Sigma$, then $\Sigma$ is Ehrhart if and only if the stellar subdivision $\Sigma_\tau$ is Ehrhart. Moreover, if both $\Sigma$ and $\Sigma_\tau$ are Ehrhart, then their Ehrhart polynomials agree: for all $f\in \PL(\Sigma)$, we have
\[
\Chi_\Sigma([f])=\Chi_{\Sigma_\tau}([f]).
\]
\end{proposition}

\begin{proof}
We first prove the ``only if'' direction by induction on dimension. To do so, assume that Ehrharticity is preserved by stellar subdivisions of fans of dimension less than $d$, and assume that $\Sigma$ is a $d$-dimensional Ehrhart fan. We aim to prove that $\Sigma_\tau$ is Ehrhart.

By Lemma~\ref{lem:starfans}, the star fan $(\Sigma_\tau)^\rho$ for any $\rho\in\Sigma_\tau(1)$ is either a star fan of $\Sigma$, a stellar subdivision of a star fan of $\Sigma$, or a fiber bundle over a star fan of $\Sigma$ with complete fibers. It follows from the Ehrharticity of $\Sigma$, the induction hypothesis, and Proposition~\ref{prop:fiberbundles} that $(\Sigma_\tau)^\rho$ is Ehrhart for all $\rho\in\Sigma_\tau(1)$. We have thus proven Condition (1) in the definition of Ehrharticity for $\Sigma_\tau$, and we turn to Condition (2).  

Given any $\tilde f\in \PL(\Sigma_\tau)$, there is a unique $f\in \PL(\Sigma)$ and $a\in\Z$ such that 
\begin{equation}
\label{eq:tildeffa}
    \tilde f=f+a\tilde\delta_{\rho_\tau}
\end{equation}
Specifically, the function $f$ takes the same value as $\tilde f$ on all ray generators except $u_{\rho_\tau}$, whereas
\[f(u_{\rho_\tau}) = \sum_{\rho \in \tau(1)} \tilde f(u_\rho)\]
to ensure linearity across $\tau$; from here, taking $a = \tilde f(u_{\rho_\tau}) - f(u_{\rho_{\tau}})$ makes \eqref{eq:tildeffa} hold.  The pictorial example in Figure~\ref{fig:ftilde} may help to clarify the relationship between $\tilde f$ and $f$.

\begin{figure}[t]
\begin{tikzpicture}[scale=1.75]
\draw[thick,fill=green!20, fill opacity=.5] (-0.75,-0.75) -- (0.75,-0.75) -- (0.75,0) -- (0,0.75) -- (-0.75,0.75) -- (-0.75,-0.75);

\draw[ultra thick, gray, opacity=0.5] (0,0) -- (0,1);
\draw[ultra thick, gray, opacity=0.5] (0,0) -- (0.9, 0.9);
\draw[ultra thick, gray, opacity=0.5] (0,0) -- (1, 0);
\draw[ultra thick, gray, opacity=0.5] (0,0) -- (0,-1);
\draw[ultra thick, gray, opacity=0.5] (0,0) -- (-1,0);

\node at (1.25,0) {$\rho_1$};
\node at (0,1.25) {$\rho_2$};
\node at (1,1) {$\rho_\tau$};

\node at (-0.75,0) {$\bullet$};
\node at (-0.5,0) {$\bullet$};
\node at (-0.25,0) {$\bullet$};
\node at (0,0) {$\bullet$};
\node at (0.25,0) {$\bullet$};
\node at (0.5,0) {$\bullet$};
\node at (0.75,0) {$\bullet$};

\node at (-0.75,0.25) {$\bullet$};
\node at (-0.5,0.25) {$\bullet$};
\node at (-0.25,0.25) {$\bullet$};
\node at (0,0.25) {$\bullet$};
\node at (0.25,0.25) {$\bullet$};
\node at (0.5,0.25) {$\bullet$};

\node at (-0.75,0.5) {$\bullet$};
\node at (-0.5,0.5) {$\bullet$};
\node at (-0.25,0.5) {$\bullet$};
\node at (0,0.5) {$\bullet$};
\node at (0.25,0.5) {$\bullet$};

\node at (-0.75,0.75) {$\bullet$};
\node at (-0.5,0.75) {$\bullet$};
\node at (-0.25,0.75) {$\bullet$};
\node at (0,0.75) {$\bullet$};

\node at (-0.75,-0.25) {$\bullet$};
\node at (-0.5,-0.25) {$\bullet$};
\node at (-0.25,-0.25) {$\bullet$};
\node at (0,-0.25) {$\bullet$};
\node at (0.25,-0.25) {$\bullet$};
\node at (0.5,-0.25) {$\bullet$};
\node at (0.75,-0.25) {$\bullet$};

\node at (-0.75,-0.5) {$\bullet$};
\node at (-0.5,-0.5) {$\bullet$};
\node at (-0.25,-0.5) {$\bullet$};
\node at (0,-0.5) {$\bullet$};
\node at (0.25,-0.5) {$\bullet$};
\node at (0.5,-0.5) {$\bullet$};
\node at (0.75,-0.5) {$\bullet$};

\node at (-0.75,-0.75) {$\bullet$};
\node at (-0.5,-0.75) {$\bullet$};
\node at (-0.25,-0.75) {$\bullet$};
\node at (0,-0.75) {$\bullet$};
\node at (0.25,-0.75) {$\bullet$};
\node at (0.5,-0.75) {$\bullet$};
\node at (0.75,-0.75) {$\bullet$};

\node at (-1,1) {$P_{\tilde{f}}$};

\end{tikzpicture}
\hspace{50bp}
\begin{tikzpicture}[scale=1.75]
\draw[thick,fill=green!20, fill opacity=.5] (-0.75,-0.75) -- (0.75,-0.75) -- (0.75,0.75) -- (-0.75,0.75) -- (-0.75,-0.75);

\draw[ultra thick, gray, opacity=0.5] (0,0) -- (0,1);
\draw[ultra thick, gray, opacity=0.5] (0,0) -- (0.9, 0.9);
\draw[ultra thick, gray, opacity=0.5] (0,0) -- (1, 0);
\draw[ultra thick, gray, opacity=0.5] (0,0) -- (0,-1);
\draw[ultra thick, gray, opacity=0.5] (0,0) -- (-1,0);

\node at (1.25,0) {$\rho_1$};
\node at (0,1.25) {$\rho_2$};
\node at (1,1) {$\rho_\tau$};

\node at (-0.75,0) {$\bullet$};
\node at (-0.5,0) {$\bullet$};
\node at (-0.25,0) {$\bullet$};
\node at (0,0) {$\bullet$};
\node at (0.25,0) {$\bullet$};
\node at (0.5,0) {$\bullet$};
\node at (0.75,0) {$\bullet$};

\node at (-0.75,0.25) {$\bullet$};
\node at (-0.5,0.25) {$\bullet$};
\node at (-0.25,0.25) {$\bullet$};
\node at (0,0.25) {$\bullet$};
\node at (0.25,0.25) {$\bullet$};
\node at (0.5,0.25) {$\bullet$};
\node at (0.75,0.25) {$\bullet$};
\draw (0.75,0.25) circle (3pt);

\node at (-0.75,0.5) {$\bullet$};
\node at (-0.5,0.5) {$\bullet$};
\node at (-0.25,0.5) {$\bullet$};
\node at (0,0.5) {$\bullet$};
\node at (0.25,0.5) {$\bullet$};
\node at (0.5,0.5) {$\bullet$};
\node at (0.75,0.5) {$\bullet$};
\draw (0.5,0.5) circle (3pt);
\draw (0.75,0.5) circle (3pt);

\node at (-0.75,0.75) {$\bullet$};
\node at (-0.5,0.75) {$\bullet$};
\node at (-0.25,0.75) {$\bullet$};
\node at (0,0.75) {$\bullet$};
\node at (0.25,0.75) {$\bullet$};
\node at (0.5,0.75) {$\bullet$};
\node at (0.75,0.75) {$\bullet$};
\draw (0.25,0.75) circle (3pt);
\draw (0.5,0.75) circle (3pt);
\draw (0.75,0.75) circle (3pt);

\node at (-0.75,-0.25) {$\bullet$};
\node at (-0.5,-0.25) {$\bullet$};
\node at (-0.25,-0.25) {$\bullet$};
\node at (0,-0.25) {$\bullet$};
\node at (0.25,-0.25) {$\bullet$};
\node at (0.5,-0.25) {$\bullet$};
\node at (0.75,-0.25) {$\bullet$};

\node at (-0.75,-0.5) {$\bullet$};
\node at (-0.5,-0.5) {$\bullet$};
\node at (-0.25,-0.5) {$\bullet$};
\node at (0,-0.5) {$\bullet$};
\node at (0.25,-0.5) {$\bullet$};
\node at (0.5,-0.5) {$\bullet$};
\node at (0.75,-0.5) {$\bullet$};

\node at (-0.75,-0.75) {$\bullet$};
\node at (-0.5,-0.75) {$\bullet$};
\node at (-0.25,-0.75) {$\bullet$};
\node at (0,-0.75) {$\bullet$};
\node at (0.25,-0.75) {$\bullet$};
\node at (0.5,-0.75) {$\bullet$};
\node at (0.75,-0.75) {$\bullet$};

\node at (-1,1) {$P_f$};

\end{tikzpicture}
\caption{A two-dimensional fan $\Sigma_\tau$ (with rays drawn in light gray) obtained from a fan $\Sigma$ by stellar subdivision along the cone $\tau$ in the first quadrant.  For a function $\tilde{f}$ taking values $\tilde{f}(u_{\rho_1}) = \tilde{f}(u_{\rho_2}) = \tilde{f}(u_{\rho_\tau}) = 3$, the polytope $P_{\tilde{f}}$ is shown on the left.  The polytope $P_f$ is shown on the right (in this case, $f(u_{\rho_{\tau}}) = 6$ and $a=-3$), and the lattice points of $P_f \setminus P_{\tilde{f}}$ are circled.}
\label{fig:ftilde}
\end{figure}
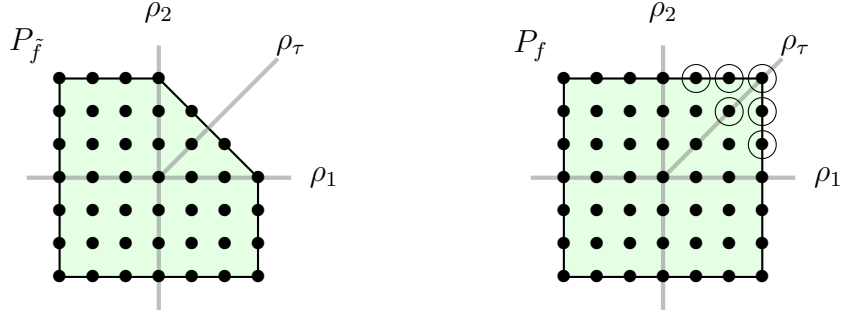

In the example of Figure~\ref{fig:ftilde}, one sees that the lattice points of $P_f$ are obtained from the lattice points of $P_{\tilde{f}}$ by adding the lattice points along the ``slices'' of $P_f \setminus P_{\tilde{f}}$ perpendicular to the ray $\rho_\tau$.  This observation motivates us to define
\[
\Chi_{\Sigma_\tau}([\tilde f])=\Chi_\Sigma([f])+\sum_{i=1}^a\Chi_{(\Sigma_\tau)^{\rho_\tau}}([f+i\tilde\delta_{\rho_\tau}]^{\rho_\tau}),
\]
where the existence of $\Chi_\Sigma$ and $\Chi_{(\Sigma_\tau)^{\rho_\tau}}$ are ensured by the assumption that $\Sigma$ is Ehrhart and the observation above that $\Chi_{(\Sigma_\tau)^{\rho_\tau}}$ is also Ehrhart.  (We adopt the standard convention that
\[\sum_{i=1}^a s_i \coloneq -\sum_{i=a+1}^0s_i\]
when $a < 1$, for any summands $s_i$; in particular in the example of Figure~\ref{fig:ftilde}, there are summands for $i=-2$, $i=-1$, and $i=0$, corresponding to the three perpendicular slices to $\rho_\tau$ along which the circled lattice points lie.)  One readily checks that $\Chi_{\Sigma_\tau}(\tilde f+m)=\Chi_{\Sigma_\tau}(\tilde f)$ for any $m\in M$, so we obtain a well-defined function $\Chi_{\Sigma_\tau}:\uPL(\Sigma_\tau)\rightarrow\Z$ that agrees with $\Chi_\Sigma$ when restricted to $\uPL(\Sigma)$.

We claim that $\Chi_{\Sigma_\tau}$ satisfies the Ehrhart relation
\begin{equation}\label{eq:blowuprelation}
\Chi_{\Sigma_\tau}([\tilde f])-\Chi_{\Sigma_\tau}([\tilde f-\tilde\delta_\rho])=\Chi_{(\Sigma_
\tau)^\rho}([\tilde f]^\rho).
\end{equation}
for every $\rho\in\Sigma_\tau(1)$. We consider four possible cases, in order of increasing complexity; for the reader's convenience, we illustrate the four cases in Figure~\ref{fig:cases}.

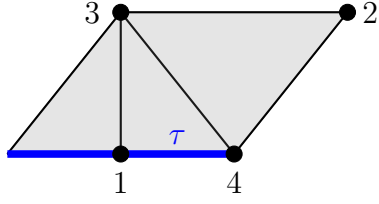
\begin{figure}[h]
\begin{tikzpicture}[scale=1.5]
\draw[thick] (0,0) -- (2,0) -- (1,1.25) -- (0,0);
\draw[thick] (1,1.25) -- (1,0);
\draw[thick] (2,0) -- (3,1.25) -- (1,1.25);

\filldraw[fill=gray, fill opacity=.2] (0,0) -- (2,0) -- (3,1.25) -- (1,1.25) -- (0,0);

\draw[blue, line width=1mm] (0,0) -- (2,0);

\node at (1.5, 0.15) {\color{blue}$\tau$};

\filldraw (1,0) circle (2pt);
\node at (1,-0.25) {$1$};

\filldraw (2,0) circle (2pt);
\node at (2,-0.25) {$4$};

\filldraw (1,1.25) circle (2pt);
\node at (0.75,1.25) {$3$};

\filldraw (3,1.25) circle (2pt);
\node at (3.2,1.25) {$2$};
\end{tikzpicture}
\caption{A 2-dimensional slice of a 3-dimensional fan $\Sigma_\tau$, with the corresponding slice of a 2-dimensional cone $\tau \in \Sigma$ highlighted.  The four labeled points represent four rays of $\Sigma_\tau$ corresponding to the four cases in the proof.}
\label{fig:cases}
\end{figure}

\begin{enumerate}
\item[Case 1:] If $\rho=\rho_\tau$, then $\tilde f=f+a\tilde\delta_{\rho_\tau}$ with $f\in \PL(\Sigma)$ implies that $\tilde f-\tilde\delta_{\rho_\tau}=f+(a-1)\tilde\delta_{\rho_\tau}$. Directly from the definition of $\Chi_{\Sigma_\tau}$, we then compute
\begin{align*}
\Chi_{\Sigma_\tau}([\tilde f])-\Chi_{\Sigma_\tau}([\tilde f-\tilde\delta_{\rho_\tau}])&=\Chi_{(\Sigma_\tau)^{\rho_\tau}}([f+a\tilde\delta_{\rho_\tau}]^{\rho_\tau})=\Chi_{(\Sigma_\tau)^{\rho_\tau}}([\tilde f]^{\rho_\tau}),
\end{align*}
verifying \eqref{eq:blowuprelation} in this case.
\item[Case 2:] If $\rho\notin\Sigma^{(\tau)}$, then we observe that
\begin{itemize}
\item $\tilde\delta_\rho=\delta_\rho$,
\item $\tilde f=f+a\tilde\delta_{\rho_\tau}$ with $f\in \PL(\Sigma)$ implies $\tilde f-\tilde\delta_\rho=(f-\delta_\rho)+a\tilde\delta_{\rho_\tau}$ with $f-\delta_\rho\in \PL(\Sigma)$, 
\item $[\tilde f]^\rho=[f]^\rho$ (since $\tilde f$ agrees with $f$ on the neighborhood of $\rho$),
\item $[f+i\tilde\delta_{\rho_\tau}]^{\rho_\tau}=[f-\delta_\rho+i\tilde\delta_{\rho_\tau}]^{\rho_\tau}$ (for a similar reason).
\end{itemize}
It then follows from the definition of $\Chi_{\Sigma_\tau}$ and the assumption that $\Sigma$ is Ehrhart that
\begin{align*}
\Chi_{\Sigma_\tau}([\tilde f])-\Chi_{\Sigma_\tau}([\tilde f-\tilde\delta_\rho])&=\Chi_{\Sigma}([f])-\Chi_{\Sigma}([f-\delta_\rho])\\
&=\Chi_{\Sigma^\rho}([f]^\rho)\\
&=\Chi_{(\Sigma_\tau)^\rho}([\tilde f]^\rho),
\end{align*}
verifying \eqref{eq:blowuprelation} in this case.
\item[Case 3:] If $\rho\in\Sigma^{(\tau)}\setminus\tau(1)$, then we observe that
\begin{itemize}
\item $\tilde\delta_\rho=\delta_\rho$, and
\item $\tilde f=f+a\tilde\delta_{\rho_\tau}$ with $f\in \PL(\Sigma)$ implies $\tilde f-\tilde\delta_\rho=(f-\delta_{\rho})+a\tilde\delta_{\rho_\tau}$ with $f-\delta_\rho\in \PL(\Sigma)$.
\end{itemize}
From the definition of $\Chi_{\Sigma_\tau}$, we then compute
\begin{align*}
\Chi_{\Sigma_\tau}([\tilde f])-\Chi_{\Sigma_\tau}([\tilde f-\tilde\delta_\rho])&=\Chi_{\Sigma}([f])-\Chi_{\Sigma}([f-\delta_\rho])\\
&\hspace{20bp}+\sum_{i=1}^a\Big(\Chi_{(\Sigma_\tau)^{\rho_\tau}}([f+i\tilde\delta_{\rho_\tau}]^{\rho_\tau})-\Chi_{(\Sigma_\tau)^{\rho_\tau}}([f-\delta_\rho+i\tilde\delta_{\rho_\tau}]^{\rho_\tau})\Big).
\end{align*}
Using the Ehrharticity of both $\Sigma$ and $(\Sigma_\tau)^{\rho_\tau}$, this expression simplifies to
\[
\Chi_{\Sigma^\rho}([f]^\rho)+\sum_{i=1}^a\Chi_{(\Sigma_\tau)^\pi}([f+i\tilde\delta_{\rho_\tau}]^\pi),
\]
where $\pi\in\Sigma_\tau(2)$ is the cone with rays $\rho$ and $\rho_\tau$. Using that $(\Sigma_\tau)^\rho=(\Sigma^\rho)_{\overline\tau}$, this final expression is then seen to be equal to $\Chi_{(\Sigma_\tau)^\rho}([\tilde f]^\rho)$, verifying \eqref{eq:blowuprelation} in this case.
\item[Case 4:] If $\rho\in\tau(1)$, then we observe that
\begin{itemize}
\item $\tilde\delta_\rho+\tilde\delta_{\rho_\tau}=\delta_\rho$ and
\item $\tilde f=f+a\tilde\delta_{\rho_\tau}$ with $f\in \PL(\Sigma)$ implies $\tilde f-\tilde\delta_\rho=(f-\delta_\rho)+(a+1)\tilde\delta_{\rho_\tau}$ with $f-\delta_\rho\in \PL(\Sigma)$.
\end{itemize}
We compute the left-hand side of \eqref{eq:blowuprelation} from the definition of $\Chi_{\Sigma_\tau}$:
\begin{align*}
\Chi_{\Sigma}([f])-\Chi_{\Sigma}([f-\delta_\rho])+\sum_{i=1}^a\Chi_{(\Sigma_\tau)^{\rho_\tau}}([f+i\tilde\delta_{\rho_\tau}]^{\rho_\tau})-\sum_{i=1}^{a+1}\Chi_{(\Sigma_\tau)^{\rho_\tau}}([f-\tilde\delta_{\rho}+(i-1)\tilde\delta_{\rho_\tau}]^{\rho_\tau}).
\end{align*}
Since $\Sigma$ is Ehrhart, the first difference above is $\Chi_{\Sigma^\rho}([f]^\rho)$, and, using the fact that the star fan $(\Sigma_\tau)^{\rho
_\tau}$ is Ehrhart, the second difference works out to be
\[
-\Chi_{(\Sigma_\tau)^{\rho_\tau}}([f-\tilde\delta_\rho]^{\rho_\tau})+\sum_{i=1}^a\Chi_{(\Sigma_\tau)^\pi}([f+i\tilde\delta_{\rho^\tau}]^\pi),
\]
where $\pi\in\Sigma_\tau(2)$ is the cone with rays $\rho$ and $\rho_\tau$. In order to verify \eqref{eq:blowuprelation} in this case, it remains to prove that
\[
\Chi_{\Sigma^\rho}([f]^\rho)-\Chi_{(\Sigma_\tau)^{\rho_\tau}}([f-\tilde\delta_\rho]^{\rho_\tau})=\Chi_{(\Sigma_\tau)^\rho}([f]^\rho).
\]
Since $(\Sigma_\tau)^\rho=(\Sigma^\rho)_{\overline\tau}$ and $[f]^\rho$ is linear on the cones of $\Sigma^\rho$, it follows that
\[
\Chi_{(\Sigma_\tau)^\rho}([f]^\rho)=\Chi_{\Sigma^\rho}([f]^\rho),
\]
so it remains to show that $\Chi_{(\Sigma_\tau)^{\rho_\tau}}([f-\tilde\delta_\rho]^{\rho_\tau})=0$. After translating $f$, we can assume that $f$ vanishes on $\tau$, so that $[f-\tilde\delta_\rho]^{\rho_\tau}=[f-\tilde\delta_\rho]$. For notational convenience, set $\Sigma'=\Pi_{\dim\tau-1}$. Using that $(\Sigma_\tau)^{\rho_\tau}$ is a $\Sigma'$-bundle over $\Sigma^\tau$, we can use the construction in the proof of Proposition~\ref{prop:fiberbundles} to write
\[
\Chi_{(\Sigma_\tau)^{\rho_\tau}}([f-\tilde\delta_\rho])=(-1)^{\dim\Sigma'}\sum_{m'\in M'}\Big(\sum_{\sigma\in\Sigma'_{-\tilde\delta_\rho',m'}}(-1)^{\dim(\sigma)}\Big)\Chi_{\Sigma^\tau}([(f-m')'']),
\]
where $(f-\tilde\delta_\rho)'=-\tilde\delta_\rho'$ follows from the assumption that $f$ vanishes on $\tau$. However, noting that the the function $-\tilde\delta_\rho'\in\PL(\Sigma')$ is concave and that the simplex $P_{\tilde\delta_\rho'}$ does not have any interior lattice points, it follows from Lemma~\ref{lem:convex} that the contribution from every $m'\in M'$ to the above sum is zero. Thus,
\[
\Chi_{(\Sigma_\tau)^{\rho_\tau}}([f-\tilde\delta_\rho]^{\rho_\tau})=0,
\]
completing the verification of \eqref{eq:blowuprelation} in this case.
\end{enumerate}

We have now proved that the stellar subdivision of an Ehrhart fan is Ehrhart, and it remains to prove the converse. By induction, assume that the Ehrharticity of a stellar subdivision implies the Ehrharticity of the original fan in all dimensions less than $d$, and let $\Sigma$ be a unimodular $d$-dimensional fan such that $\Sigma_\tau$ is Ehrhart for some $\tau\in\Sigma$. We must prove that $\Sigma$ is Ehrhart.

To verify the first Ehrhart condition for $\Sigma$, we note that, for any $\rho\in\Sigma$, the star fan $\Sigma^\rho$ is either equal to a star fan of $\Sigma_\tau$ or has a stellar subdivision that is equal to a star fan of $\Sigma_\tau$. Using the assumption that $\Sigma_\tau$ is Ehrhart and the induction hypothesis, it follows that $\Sigma^\rho$ is Ehrhart for any $\rho\in\Sigma(1)$.

To verify the second condition of Ehrharticity for $\Sigma$, start by defining $\Chi_\Sigma:\uPL(\Sigma)\rightarrow\Z$ by
\[
\Chi_\Sigma([f])\coloneq \Chi_{\Sigma_\tau}([f]),
\]
where we are using the assumption that $\Sigma_\tau$ is Ehrhart. Thus, Condition (2) of Ehrharticity is equivalent to proving that
\begin{equation}\label{eq:blowdownrelation}
\Chi_{\Sigma_\tau}([f])-\Chi_{\Sigma_\tau}([f-\delta_\rho])=\Chi_{(\Sigma_
\tau)^\rho}([f]^\rho)
\end{equation}
for every $f\in \PL(\Sigma)$ and $\rho\in\Sigma(1)$. We consider two cases.

\begin{enumerate}
\item[Case 1:] If $\rho\notin\tau(1)$, then $\delta_\rho=\tilde\delta_\rho$, and it follows from the Ehrharticity of $\Sigma_\tau$ that
\[
\Chi_{\Sigma_\tau}([f])-\Chi_{\Sigma_\tau}([f-\delta_\rho])=\Chi_{\Sigma_\tau}([f])-\Chi_{\Sigma_\tau}([f-\tilde\delta_\rho])=\Chi_{(\Sigma_
\tau)^\rho}([f]^\rho)
\]
verifying \eqref{eq:blowdownrelation} in this case.
\item[Case 2:] If $\rho\in\tau(1)$, then we observe that $\delta_\rho=\tilde\delta_\rho+\tilde\delta_{\rho_\tau}$ and it follows that
\[
\Chi_{\Sigma_\tau}([f])-\Chi_{\Sigma_\tau}([f-\delta_\rho])=\Chi_{\Sigma_\tau}([f])-\Chi_{\Sigma_\tau}([f-\tilde\delta_\rho-\tilde\delta_{\rho_\tau}]).
\]
Ehrharticity of $\Sigma_\tau$ implies that
\[
\Chi_{\Sigma_\tau}([f])-\Chi_{\Sigma_\tau}([f-\tilde\delta_\rho-\tilde\delta_{\rho_\tau}])=\Chi_{(\Sigma_\tau)^{\rho}}([f]^\rho)+\Chi_{(\Sigma_\tau)^{\rho_\tau}}([f-\tilde\delta_\rho]^{\rho_\tau}).
\]
The argument used in Case 4 above shows that 
\[
\Chi_{(\Sigma_\tau)^{\rho_\tau}}([f-\tilde\delta_\rho]^{\rho_\tau})=0,
\]
verifying \eqref{eq:blowdownrelation} in this case, and completing the proof of the proposition.\qedhere
\end{enumerate}
\end{proof}

\subsection{Ehrharticity is independent of fan structure}

We now come to the main theorem of this section, which states that Ehrharticity and the values of the Ehrhart polynomial only depend on the support of a fan.

\begin{theorem}\label{thm:support}
Let $\Sigma_1$ and $\Sigma_2$ be unimodular fans with $|\Sigma_1|=|\Sigma_2|$. Then $\Sigma_1$ is Ehrhart if and only if $\Sigma_2$ is Ehrhart. Moreover, if $\Sigma_1$ and $\Sigma_2$ are Ehrhart and $f\in \PL(\Sigma_1)\cap \PL(\Sigma_2)$, then
\[
\Chi_{\Sigma_1}([f])=\Chi_{\Sigma_2}([f]).
\]
\end{theorem}

The proof of Theorem~\ref{thm:support} utilizes the following important result of W{\l}odarczyk, often referred to as the \emph{Weak Factorization Theorem}.

\begin{theorem}[\cite{wlodarczyk} Theorem~13.3]\label{thm:wlodarczyk}
Let $\Sigma$ be any simplicial fan and let $\Sigma'$ and $\Sigma''$ be two unimodular refinements of $\Sigma$. There exists a sequence $\Sigma_0,\dots,\Sigma_k$ of unimodular refinements of $\Sigma$ such that $\Sigma_0=\Sigma'$, $\Sigma_k=\Sigma''$, and for each $i=0,\dots,k-1$, either $\Sigma_i$ is a stellar subdivision of $\Sigma_{i+1}$ or $\Sigma_{i+1}$ is a stellar subdivision of $\Sigma_i$.
\end{theorem}

\begin{proof}[Proof of Theorem~\ref{thm:support}]
Let $\Sigma$ be a common unimodular refinement of $\Sigma_1$ and $\Sigma_2$. Theorem~\ref{thm:wlodarczyk} implies that there exists a sequence of unimodular fans $\Sigma_1',\dots,\Sigma_k'$ such that
\[
\Sigma_1\sim\Sigma_1'\sim\dots\sim\Sigma_j'\sim\Sigma\sim\Sigma_{j+1}'\sim\dots\sim\Sigma_k'\sim\Sigma_2,
\]
where each $\sim$ denotes that one fan is a stellar subdivision of the other---but is agnostic about which direction the stellar subdivision goes---and where $\Sigma_1',\dots,\Sigma_j'$ all refine $\Sigma_1$ and $\Sigma_{j+1}',\dots,\Sigma_k'$ all refine $\Sigma_2$. It then follows from Proposition~\ref{prop:blowups} that any fan in this sequence is Ehrhart if and only if any other fan in this sequence is Ehrhart; in particular, $\Sigma_1$ is Ehrhart if and only if $\Sigma_2$ is Ehrhart. Morever, if any fan in this sequence is Ehrhart and if $f\in\PL(\Sigma_1)\cap \PL(\Sigma_2)$, then $f\in\PL(\Sigma_i')$ for all $i=1,\dots,k$, and Proposition~\ref{prop:blowups} further implies that the Ehrhart polynomial of any fan in this sequence takes the same value at $f$; in particular, $\Chi_{\Sigma_1}([f])=\Chi_{\Sigma_2}([f])$.
\end{proof}

Theorem~\ref{thm:support} allows us to make the following definition.

\begin{definition}\label{def:Ehrhartfansupport}
Let $X$ be the support of an Ehrhart fan. The \textbf{Ehrhart polynomial of $X$} is the function $\Chi_X:\uPL(X)\rightarrow\Z$ defined by
\[
\Chi_X([f])=\Chi_\Sigma([f])
\]
where $\Sigma$ is any unimodular fan with support $X$.
\end{definition}

\subsection{Ehrhart polynomials and piecewise-exponential functions}

We now show that Ehrhart polynomials induce well-defined linear functions on rings of piecewise-exponential functions, generalizing the fact that lattice-point counts induce a well-defined linear function on the polytope algebra of lattice polytopes. 

\begin{theorem}
\label{thm:chiK}
If $X$ is the support of an Ehrhart fan, then the Ehrhart polynomial of $X$ induces a well-defined linear function $\widetilde\Chi_X:\uPE(X)\rightarrow \Z$ defined on additive generators by
\[
\widetilde\Chi_X([e^{f}]) = \Chi_X([f]).
\]
\end{theorem}

\begin{proof}
By Theorem~\ref{thm:K(X)presentation}, it suffices to prove that
\begin{equation}\label{eq:rel1}
\Chi_X([f])+\Chi_X([g])=\Chi_X([\max(f,g)])+\Chi_X([\min(f,g)])
\end{equation}
for all $f,g\in\PL(X)$, and that
\begin{equation*}
\Chi_X([f+\ell])=\Chi_X([f])
\end{equation*}
for all $f\in\PL(X)$ and $\ell\in\L(X)$, the latter of which is immediate from the observation that $[f+\ell]=[f]\in\PL(X)$. To prove \eqref{eq:rel1}, start by fixing a unimodular fan $\Sigma$ supported on $X$ such that $f$, $g$, $\max(f,g)$, and $\min(f,g)$ are all linear on the cones of $\Sigma$. Notice that, for each cone $\sigma\in\Sigma$, it must be the case that $f\geq g$ or $f\leq g$ on all of $\sigma$. By Definition~\ref{def:Ehrhartfansupport}, we see that \eqref{eq:rel1} is equivalent to the similar relation obtained by replacing $X$ with $\Sigma$ everywhere:
\begin{equation}\label{eq:rel3}
\Chi_\Sigma([f])+\Chi_\Sigma([g])=\Chi_\Sigma([\max(f,g)])+\Chi_\Sigma([\min(f,g)]).
\end{equation}
We prove \eqref{eq:rel3} by induction on the nonnegative integer
\begin{equation}\label{eq:difference}
\sum_{\rho\in\Sigma(1)}|f(u_\rho)-g(u_\rho)|.
\end{equation}

When \eqref{eq:difference} is zero, then $f=g=\max(f,g)=\min(f,g)$, and \eqref{eq:rel3} follows immediately, proving the base case. To prove the induction step, suppose that \eqref{eq:difference} is positive and choose any $\rho\in\Sigma(1)$ such that (without loss of generality) $f(u_\rho)>g(u_\rho)$. By our assumption on $\Sigma$, it follows that $f(x)\geq g(x)$ for all $x$ in the neighborhood of $\rho$. In particular, it then follows that
\[
[f]^\rho=[\max(f,g)]^\rho\in\uPL(\Sigma^\rho).
\]
and
\[
\max(f-\delta_\rho,g)=\max(f,g)-\delta_\rho.
\]
Applying the Ehrhart recursion to both $[f]$ and $[\max(f,g)]$ at $\rho$, we then obtain
\[
\Chi_\Sigma([f])=\Chi_\Sigma([f-\delta_\rho])+\Chi_{\Sigma^\rho}([f]^\rho)
\]
and
\[
\Chi_\Sigma([\max(f,g)])=\Chi_\Sigma([\max(f-\delta_\rho,g)])+\Chi_{\Sigma^\rho}([f]^\rho).
\]
In particular, it then follows that \eqref{eq:rel3} is equivalent to 
\[
\Chi_\Sigma([f-\delta_\rho])+\Chi_\Sigma([g])=\Chi_\Sigma([\max(f-\delta_\rho,g)])+\Chi_\Sigma([\min(f,g)]),
\]
but this new relation has a strictly smaller value of \eqref{eq:difference}, and thus follows by the induction hypothesis.
\end{proof}

Upon precomposing the linear function of Theorem~\ref{thm:chiK} with the natural ring homomorphism $\phi:\uPE(\Sigma)\rightarrow\uPE(|\Sigma|)$, we obtain the following consequence.

\begin{corollary}
If $\Sigma$ is an Ehrhart fan, then the Ehrhart polynomial of $\Sigma$ induces a well-defined linear function $\widetilde\Chi_\Sigma:\uPE(\Sigma)\rightarrow \Z$ defined on additive generators by
\[
\widetilde\Chi_\Sigma([e^{f}]) = \Chi_\Sigma([f]).
\]
\end{corollary}

\section{Bergman fans are Ehrhart}
\label{sec:matroidfans}

In this final section of the paper, we prove that the Bergman fan $\Sigma_\M$ of a (loopless) matroid $\M$ is Ehrhart, where the Ehrhart polynomial is derived from the ``Euler characteristic'' of the matroid. This completes our goal of exhibiting Ehrhart functions on Ehrhart fans as a unifying framework for Euler characteristics of both smooth complete toric varieties and matroids. We begin by recalling the relevant notions from the theory of matroids and their Euler characteristics.

\subsection{Matroids and Bergman fans}

A matroid, by definition, consists of a finite set $E$ called the {\bf ground set}, together with a distinguished collection of subsets of $E$ called the {\bf independent sets}, which are required to satisfy certain axioms.  The prototypical example is the case where $E$ is a finite set of vectors in some vector space and the independent sets are precisely those subsets of $E$ that are linearly independent.  Motivated by this example, when given a matroid $\M$ on ground set $E$, one can define the {\bf rank} of any subset $A \subseteq E$ as the maximum size of an independent subset of $A$.  A subset $F \subseteq E$ is said to be a {\bf flat} of $\M$ if adding any additional element of $E$ to $F$ increases its rank; in the prototypical example mentioned above, flats are the intersections of $E$ with linear subspaces of the ambient vector space. Note that the full ground set $E$ is a flat, and we say that $\M$ is \textbf{loopless} if the empty set is also a flat. (The terminology comes from another fundamental example of a matroid, in which the ground set is the edge set of a fixed finite graph and a subset is declared to be independent if it contains no cycles.) We use the terminology {\bf proper flat} to refer to any flat of a loopless matroid $\M$ that is neither empty nor all of $E$.

In what follows, we will denote the ground set by
\[
E = \{0,1,2,\ldots, n\}
\]
for some integer $n\ge 0$.  Let $\M$ be a matroid on this ground set, and let $N = \Z^E/\Z(1,\ldots,1).$ The {\bf Bergman fan} $\Sigma_\M$ is a rational fan in the vector space
\[
N_\R := \R^E/\R(1,\ldots,1) \cong \R^n.
\]
To define its cones, for each $i \in E$, let $\overline{\e}_i \in N_\R$ denote the image under the above quotient of the $i$th standard basis vector in $\R^E$. Then, for each proper flat $F$ of $\M$, let
\[
\overline{\e}_F := \sum_{i \in F} \overline{\e}_i \in N_\R.
\]
There is one cone $\sigma_{\mathcal{F}}$ of $\Sigma_\M$ for each chain $\mathcal{F} = (F_1 \subsetneq F_2 \subsetneq \cdots \subsetneq F_k)$ of proper flats of $\M$, given by
\[
\sigma_{\mathcal{F}} = \text{Cone}(\overline{\e}_{F_1}, \overline{\e}_{F_2},\ldots, \overline{\e}_{F_k}). 
\]
For each proper flat $F$, we denote the associated ray of $\Sigma_\M$ by $\rho_F$. It is not hard to verify from the definitions that the Bergman fan $\Sigma_\M$ is unimodular with respect to $N$.

\subsection{$\K$-rings and Euler characteristics of matroids}
\label{subsec:K(M)}

The unimodular fan $\Sigma_\M$ has an associated smooth toric variety $Z_{\Sigma_\M}$, and Larson--Li--Payne--Proudfoot define the \textbf{matroid $\K$-ring} to be the $\K$-ring (that is, the Grothendieck ring of vector bundles) of this toric variety:
\[
K(\M)\coloneq K(Z_{\Sigma_\M}).
\]
Larson--Li--Payne--Proudfoot also give a concrete presentation for $\K(\M)$ as a quotient of the Stanley--Reisner ring of $\Sigma_\M$, but that presentation is not necessary for what follows. What is important to note is that $K(\M)$ is generated as a ring by the line bundles $\mathcal{O}(D_F)$, where $D_F$ denotes the toric divisor associated to the ray $\rho_F\in\Sigma_\M$; these generators are referred to in \cite{LLPP} as the ``Feichtner--Yuzvinsky" generators of $\K(\M)$. Motivated by Euler characteristics of wonderful compactifications in the representable case, Larson--Li--Payne--Proudfoot then define the \textbf{matroid Euler characteristic} for any loopless matroid $\M$, which is an explicit linear map
\[
\Chi_\M: K(\M)\rightarrow\Z,
\]
and which agrees with the usual Euler characteristic of vector bundles on the associated wonderful compactification whenever $\M$ is representable. The relevant properties of $\Chi_\M$ that we require are cited in the proof of Theorem~\ref{thm:matroid} below. 

\subsection{Bergman fans are Ehrhart}

We now come to the main result of this section.

\begin{theorem}
\label{thm:matroid}
Let $\M$ be a loopless matroid.  Then $\Sigma_\M$ is Ehrhart, with Ehrhart polynomial given by
\[
\Chi([f]) = \Chi_\M\left(\bigotimes_{F}\O(D_F)^{\otimes f(u_F)}\right).
\]
\end{theorem}
\begin{proof}
We first note that $\Chi$ is well-defined on $\uPL(\Sigma_\M)$ simply because 
\[
\uPL(\Sigma_\M)\cong \mathrm{Pic}(Z_{\Sigma_\M})\subseteq K(\Sigma_{\M})= K(\M),
\]
where the isomorphism is given by (\cite[Theorems~4.2.1 and 4.2.12]{CoxLittleSchenck})
\[
[f]\mapsto \bigotimes_{F}\O(D_F)^{\otimes f(u_F)}.
\]

Next, for any ray $\rho_F$ of $\Sigma_\M$ corresponding to a proper flat $F$ of $\M$, the star fan of $\Sigma_\M$ at $\rho_F$ is isomorphic to a product of Bergman fans,
\begin{equation}
    \label{eq:Bergmanstarfan}
    \Sigma_\M^{\rho_F} \cong \Sigma_{\M^F} \times \Sigma_{\M_F},
\end{equation}
by \cite[Proposition 3.5]{AHK}.  Here, $\M^F$ denotes the restriction of $\M$ to $F$ and $\M_F$ denotes the contraction, both of which are themselves loopless matroids. Thus, condition (1) in the definition of Ehrharticity follows by induction on the size of the ground set and the fact that products of Ehrhart fans are Ehrhart (Proposition~\ref{prop:products}), where the base case is given by the (zero-dimensional) Bergman fan of the unique loopless matroid on a one-element ground set.

The fact that $\Chi$ satisfies the recursion in condition (2) of the definition of Ehrharticity follows directly from \cite[Proposition 8.6]{LLPP} via the identification \eqref{eq:Bergmanstarfan}. To help the reader parse notation, we note that \cite{LLPP} uses the notation 
$\mathrm{Snap}_\M^{\mathrm{FY}}(\mathbf{a})$ for a tuple $\mathbf{a} = (a_F)$ of nonnegative integers, which is related to the above by
\[
\mathrm{Snap}_\M^{\mathrm{FY}}(\mathbf{a})=\Chi_\M\left(\bigotimes_{F}\O(D_F)^{\otimes a_F}\right).
\]
And finally, that $\Chi(0)=1$ follows from setting $\mathbf{m}=0$ in \cite[Theorem~7.2]{LLPP}.
\end{proof}

\bibliographystyle{alpha}
\bibliography{references}

\end{document}